\documentclass[10.5pt]{article}
\usepackage[leqno]{amsmath}

\usepackage{graphicx}

\usepackage{amsmath}
\usepackage{amssymb}
\usepackage{latexsym}
\usepackage{amsmath, amsfonts,amssymb, amsthm, euscript,makeidx,color,mathrsfs}

\usepackage{cases}

\usepackage{bm}
\usepackage{latexsym,enumerate,amscd,verbatim}

\def\sqr#1#2{{\vcenter{\vbox{\hrule height.#2pt
              \hbox{\vrule width.#2pt height#1pt \kern#1pt \vrule width.#2pt}
              \hrule height.#2pt}}}}

\def\3n{\negthinspace \negthinspace \negthinspace }
\def\2n{\negthinspace \negthinspace }
\def\1n{\negthinspace }

\def\dbB{\mathbb{B}}

\def\dbE{\mathbb{E}}
\def\dbF{\mathbb{F}}

\def\dbH{\mathbb{H}}

\def\dbN{\mathbb{N}}

\def\dbP{\mathbb{P}}
\def\dbQ{\mathbb{Q}}
\def\dbR{\mathbb{R}}

\def\dbX{\mathbb{X}}

\def\sB{\mathscr{B}}

\def\sE{\mathscr{E}}

\def\sH{\mathscr{H}}

\def\sL{\mathscr{L}}

\def\sP{\mathscr{P}}

\def\sX{\mathscr{X}}

\def\={\buildrel \triangle \over =}

\def\ds{\displaystyle}

\def\ns{\noalign{\ss}}
\def\a{\alpha}
\def\b{\beta}

\def\d{\delta}
\def\e{\varepsilon}

\def\l{\lambda}

\def\si{\sigma}
\def\t{\tau}
\def\f{\varphi}

\def\o{\omega}

\def\D{\Delta}
\def\Th{\Theta}
\def\L{\Lambda}

\def\cA{{\cal A}}
\def\cB{{\cal B}}

\def\cF{{\cal F}}
\def\cG{{\cal G}}
\def\cH{{\cal H}}

\def\cL{{\cal L}}
\def\cM{{\cal M}}
\def\cN{{\cal N}}

\def\cQ{{\cal Q}}

\def\cU{{\cal U}}

\def\cX{{\cal X}}

\def\ss{\smallskip}

\def\mds{\medskip}

\def\q{\quad}
\def\qq{\qquad}
\def\hb{\hbox}

\def\lan{\mathop{\langle}}
\def\ran{\mathop{\rangle}}

\def\h{\widehat}
\def\wt{\widetilde}

\def\cd{\cdot}

\def\({\Big (}
\def\){\Big )}
\def\[{\Big[}
\def\]{\Big]}
\def\bde{\begin{definition}}
\def\ede{\end{definition}}
\def\be{\begin{equation}}
\def\bel{\begin{equation}\label}
\def\ee{\end{equation}}
\def\bt{\begin{theorem}}
\def\et{\end{theorem}}
\def\bc{\begin{corollary}}
\def\ec{\end{corollary}}
\def\bl{\begin{lemma}}
\def\el{\end{lemma}}
\def\bp{\begin{proposition}}
\def\ep{\end{proposition}}
\def\bas{\begin{assumption}}
\def\eas{\end{assumption}}
\def\br{\begin{remark}}
\def\er{\end{remark}}
\def\ba{\begin{array}}
\def\ea{\end{array}}
\def\ed{\end{document}}

\def\square#1{\vbox{\hrule\hbox{\vrule height#1%
     \kern#1\vrule}\hrule}}
\def\rectangle#1#2{\vbox{\hrule\hbox{\vrule height#1%
     \kern#2\vrule}\hrule}}

\font\tenbb=msbm10 \font\sevenbb=msbm7 \font\fivebb=msbm5

\newfam\bbfam
\scriptscriptfont\bbfam=\fivebb \textfont\bbfam=\tenbb
\scriptfont\bbfam=\sevenbb

\newtheorem{lemma}{Lemma}[section]
\newtheorem{remark}{Remark}[section]

\newtheorem{theorem}{Theorem}[section]
\newtheorem{corollary}{Corollary}[section]

\newtheorem{definition}{Definition}[section]
\newtheorem{proposition}{Proposition}[section]
\newtheorem{assumption}{Assumption}[section]

\makeatletter
   
   \@addtoreset{equation}{section}
\makeatother

\makeatletter

\def\mmf{\mathbb{F}}

\def\dbQ{\mathbb{Q}}

\def\t{\tau}

\renewcommand{\@seccntformat}[1]{\csname the#1\endcsname.\hspace{0.5em}}
\makeatother

\author{Tianxiao Wang\footnote{School of Mathematics, Sichuan University,
Chengdu, Sichuan Province, 610065, China. Email: wtxiao2014@scu.edu.cn.
The research was supported by NSF of China under grant
11231007, 11301298, 11401404 and 11471231, China Postdoctoral Science Foundation (2014M562321).}}

\begin{document}
\title{General maximum principles for optimal control problems of stochastic Volterra integral equations}

\date{}
\maketitle

\begin{abstract}
Optimal control problems of forward stochastic Volterra
integral equations (SVIEs) are formulated and studied. When control region is arbitrary subset of Euclidean space and control enters into the diffusion, necessary conditions of Pontryagin's type for optimal controls are established via spike variation. Our conclusions naturally cover the analogue of stochastic differential equations (SDEs), and our developed methodology drops the reliance on It\^{o} formula and second-order adjoint equations.
Some new features, that are concealed in the SDEs framework, are revealed in our situation. For example, instead of using second-order adjoint equations, it is more
appropriate to introduce $second$-$order$ $adjoint$ $processes$.
Moreover, the conventional way of using $one$ second-order adjoint equation is inadequate here. In other words, $two$ adjoint processes, which just merge into the solution of second-order adjoint equation in SDEs situation, are actually required and proposed in our setting.

\end{abstract}

\vspace{13mm}


\bigskip


\maketitle

\section{Introduction}

Suppose $(\Omega,\cF, \mmf,$ $P)$ is a complete probability space, $W(\cdot)$ is a
one-dimensional Wiener process which generates filtration $\mmf=\{\cF_{t} \}_{0\le t\le T}$.
In this paper, we study optimal control problems for stochastic Volterra integral equations (SVIEs, in short), where the state equation is described as,
\bel{FSVIE-3-1}\2n\1n\ba{ll}
\ns\ds
X(t)=\f(t)+\int_0^tb(t,s,X(s),u(s))ds+\int_0^t\si(t,s,X(s),u(s))dW(s),
\ea\ee
with $t\in[0,T]$ and the cost functional is
\bel{cost3.1}
J(u(\cd))= \dbE\[h(X(T)+\int_0^Tl(s,X(s),u(s))ds\].
\ee
Here $u(\cd)$ is a $control$ $process$ taking values in $U\subset\dbR^m$, and $X(\cd)$, the (strong) solution of (\ref{FSVIE-3-1}), is the corresponding $state$ $process$ in $\dbR^n$. Under proper conditions, (\ref{FSVIE-3-1}) admits a unique solution $X(\cd)$. Then the optimal control problem is to find suitable $u(\cd)$ to minimize (\ref{cost3.1}).

One usual way to treat above problem is to establish the Pontryagin's maximum principle. In 1964, Friedman \cite{Friedman 1964} discussed the case of
\bel{deterministic-VIE-intro}\2n\1n\ba{ll}
\ns\ds
X(t)=x_0+\int_0^th(t-s)b(s,X(s),u(s))ds,\qq t\in[0,T].
\ea\ee
Later in 1969, Vinokurov \cite{Vinokurov 1969} investigated the case of general nonlinear VIEs with constrained state processes. Some relevant works in the deterministic framework can be found in e.g. Bakke \cite{Bakke 1974}, Bonnans et al \cite{Bonnans et al 2013}, Dmitruk-Osmolovski \cite{Dmitruk 2014}, Halanay \cite{Halanay 1968} and the references therein. As to the stochastic case, Yong \cite{Yong 2006}, \cite{Yong 2008} derived the maximum principle of optimal control problems when $U$ is convex. {\O}ksendal-Zhang \cite{Oksendal-Zhang 2010}, Agram-{\O}ksendal \cite{Agram-Oksendal 2015} present some investigations with open set $U$ by means of Malliavin calculus. Some related studies include Bonaccorsi et al \cite{Bonaccorsi-Confortola-Mastrogiacomo 2012}, Shi et al \cite{Shi-Wang-Yong 2015}, Wang-Zhang \cite{Wang-Zhang 2016} and so on. However, none of above papers can treat the case when both $U:=\{0,1\}$ and diffusion depends on control. This aim of this current paper is to spread out detailed discussions with arbitrary $U\subset\dbR^m$ and control-dependent $\si $.

To explain the motivations of studying SVIE (\ref{FSVIE-3-1}), we start with some comparisons with classical stochastic differential equations (SDEs, in short). If $\f,\ b,\ \si$ are independent of $t$, then for $t\in[0,T]$, $X(t)$ satisfies the following controlled SDE
\bel{FSDE-state-equation-intro} \!\!\!\!\!\!\!\!  \ba{ll}
\ns\ds
 X(t)\!\!=x_0+\int_0^tb(s,X(s),u(s))ds+\int_0^t\si(s,X(s),u(s))dW(s).
\ea\ee
We list out three new features of (\ref{FSVIE-3-1}) which can not captured by above (\ref{FSDE-state-equation-intro}).

$\lozenge$ The diffusion could be $\cF_t$-measurable, and thus involves anticipated stochastic integral. Even so, one can still obtain adapted solutions, see e.g. \cite{Pardoux-Protter 1990}.

$\lozenge$  Due to the reliance on both $t$ and $s$, the drift or diffusion could have singular kernels, which might bring new conclusions and challenges, see e.g. \cite{Cochran et al 1995}, \cite{De-Hoog-Weiss 1973}, \cite{Lin-Yong-2018}.

$\lozenge$  (\ref{FSDE-state-equation-intro}) is memoryless in the sense that $X(t+\D t)-X(t)$ only depends on the values of $b,\ \si$ in $[t,t+\D t]$. Nevertheless, in reality, long-term dependence usually exists, and (\ref{FSVIE-3-1}) is one proper choice to represent the memory effect, see e.g. \cite{Agram-Oksendal 2015}, \cite{Shi-Wang-Yong 2015}, \cite{Wang-2016-submitted}.

Based on these facts, we believe that optimal control problems of (\ref{FSVIE-3-1}) are much richer than that of (\ref{FSDE-state-equation-intro}). Besides the theoretical parts, we emphasize that stochastic (or deterministic) VIEs can describe many specific models, such as optimal dynamic advertising model (p.53-p.55 in Hartl \cite{Hartl 1984}), optimal capital policy model (p.23-p.25 of Arrow \cite{K. L Arrow 1964} or p.469, p.472 of Kamien-Muller \cite{Kamien-Muller 1976}), Ramsey vintage capital model (p.582 of Hritonenko-Yatsenko \cite{Hritonenko-Yatsenko 2008}), stochastic subdiffusion phenomenon in biophysics experiments (p.506-p.507 of Kou \cite{S.C. Kou 2008}), stochastic inventory-production model (p.2576-p.2578 of \cite{Wang-Zhang 2016}), optimal investment model with memory (p.1088-p.1093 of \cite{Agram-Oksendal 2015}), capital stock model in economics (Example 3.1 of Pardoux-Protter \cite{Pardoux-Protter 1990}).

Next we introduce one more stochastic epidemic prevention model.
Suppose there is a population suffered from one infected disease during $[0,T]$, where $T$ is the time when there is no infective individual, and $0$ is the time when the infective group are separated and received medical treatment.
At time $t\in[0,T]$, we denote by $X(t)$ the population of infected people (including the ones who received vaccines before recovering), $u(t)$ the amount of vaccines provided by local government. Suppose any infected individual is likely to get worse and die at future time. Hence we define
random variable $\xi_1 $ his/her life length, $m_1(r)dr$ his/her dying probability during $[r,r+dr]$ with density function $m_1(\cd)$. Consequently,
$$ F_1(r):=\dbP\big\{\xi_1\leq r\big\}= \int_0^{r}m_1(s)ds, \ \ r\in[0,T].
$$
For $t\in[0,T]$, $s\in[0,t]$, at time $t-s$ the infective population is $X(t-s)$. After time $s$ (i.e. at future time $t$) the dying population becomes $X(t-s)m_1(s)ds$.
Then, the total number between $[0,t]$ is
$$\int_0^tX(t-s)m_1(s)ds=\int_0^t X(s)m_1(t-s)ds.$$
We shall observe that the infected people in different stage (or time period) respond to the vaccine in distinctive ways. Therefore, it is reasonable to introduce the efficiency index $a(\cd)$ depending on time. As a result, at time $t$ there are $a(t)u(t)$ amount of people whose scenarios become stable.
For this group, we define random variables $\xi_2$ the duration of recovering completely, $m_2(r)dr$ the probability of getting normal during $[r,r+dr]$ with density function $m_2(\cd)$. In other words,
$$F_2(r):= \dbP\big\{\xi_2\leq r\big\}= \int_0^{r}m_2(s)ds,\ \ r\in[0,T].
$$
At time $t-s$ with $s\in[0,t]$, the amount of vaccines is $u(t-s)$, and the population with stable physical condition is $a(t-s)u(t-s)$. After time $s$, there are $a(t-s)u(t-s)m_2(s)ds$ amount of people healing from the disease. Thus the total number between $[0,t]$ is
$$\int_0^ta(t-s)u(t-s)m_2(s)ds=\int_0^ta(s) u(s)m_2(t-s)ds.$$
To sum up, the increment of infected individuals at time $t$ is
$$\ba{ll}
\ns\ds \D X(t)=\Big[-\int_0^tm_1(t-s)X(s)ds-\int_0^tm_2(t-s)a(s)u(s)ds\Big]\D t.
\ea$$
In other words, for $t\in[0,T]$,
$$\ba{ll}
\ns\ds X(t)=x_0 -\int_0^t \Big[\int_0^sm_1(s-r)X(r)dr\Big]ds-\int_0^t \Big[\int_0^sm_2(s-r)a(r)u(r)dr\Big]ds\\
\ns\ds\qq=x_0-\int_0^t F_1(s-r)X(r)dr-\int_0^t F_2(s-r)a(r)u(r)dr ,
\ea$$
where $x_0$ is the infected individuals at time $0$. At this very moment, we make some points on the efficiency coefficient $a(\cd)$. Observe that it can be easily influenced by other random factors, like the individual's physical quality, the epidemic situation in global group, the improvement of the vaccine, etc.  Consequently, we may replace $a(t)$ with $a(t)+\dot{W}(t)$, where $\dot{W}(\cd)$ represents the white noise. Then for $t\in[0,T]$
\bel{Linear-SVIE-in-Example}\ba{ll}
\ns\ds X(t)=x_0\!\!-\!\!\int_0^t\!\!\big[F_1(t-r)X(r)\!+\!F_2(t-r)a(r) u(r)\big]dr \!\!-\!\!\int_0^tF_2(t-r)u(r)dW(r).
\ea\ee
Suppose the local government wants to find suitable $\bar u(\cd)$ to minimize
\bel{Quadratic-cost-in-Example}\ba{ll}
\ns\ds J(u(\cd)):=\dbE \int_0^T \big(G_1(X(s))+G_2(u(s))\big)ds,
\ea\ee
where $G_1(\cd)$ is the daily cost of living for the infected people, $\ G_2(\cd)$ represents the research and development cost with respect to vaccines. Hence we come up with an optimal control problem associated with (\ref{Linear-SVIE-in-Example}), (\ref{Quadratic-cost-in-Example}).

We return back to our optimal control problem associated with (\ref{FSVIE-3-1}), (\ref{cost3.1}). When (\ref{FSVIE-3-1}) reduces to (\ref{FSDE-state-equation-intro}), $U\subset\dbR^{m}$ is arbitrary and $\sigma$ depends on $u(\cd)$, the maximum principles of optimal controls were firstly solved in Peng \cite{Peng90} with spike variation.
We emphasize that our investigations here are by no means fairly straightforward generalization, and more fresh thoughts have to be injected due to the encountered challenges. We refer to Subsection 3.1 and Section 4 for more details.

The rest of this paper is organized as follows. In Section 2, some notations, spaces are introduced and the optimal control problem is formulated in detail. Section 3 includes four parts. The first part is aim to illustrate how the obstacles of our study arise, as well as some intuitive introductions of our developed approach. The second part and third part are devoted to treating the encountered difficulties. In the fourth part, two main results of this paper are established and several special cases are discussed. In Section 4, some concluding remarks are present. Finally, a few key lemmas are given in the Appendix.

\section{Preliminaries}

\subsection{Some notations}

First of all, let us introduce some spaces. For
$H:=\dbR,\dbR^n,\dbR^{n\times m}$, etc., we denote its norm by $|\cd|$.
For $0\le s<t\le T$, $p>1$, we define
$$\ba{ll}
\ns\ds L^p_{\cF_t}(\Omega;H):=\Big\{X:\Omega\to H\bigm|X\hb{ is
$\cF_t$-measurable, }\dbE|X|^p<\infty\Big\},\\
%
\ns\ds C_{\dbF}\big([s,t];L^p(\Omega;H)\big)\1n:=\1n\Big\{X\1n:\1n[s,t]\times\Omega\to
H\bigm|X(\cd)\hb{ is
continuous from $[s,t]$ to}\\
\ns\ds\qq\qq\qq\qq\qq\ L^2(\Omega;H), \hb{and measurable, $\dbF$-adpated},\  \sup_{r\in[s,t]}\dbE|X(r)|^p<\infty\Big\},\\
\ns\ds L^p_\dbF(s,t;H)\!\!:=\!\!\Big\{\!X\!\!:\1n[s,t]\times\Omega\!\!\to
H \!\!\bigm|\!\!X(\cd)\hb{ is $\dbF$-adapted, measurable}, \dbE\!\!\int_s^t\!\!|X(r)|^pdr\!\!<\!\!\infty\Big\}, \\
\ns\ds L^2(s,t;L^2_\dbF(s,t;H)):=\!\Big\{Z:[s,t]^2\times\Omega\!\to\!
H\!\bigm|\! Z(\cd,\cd)\ \hb{is measurable,} \ \hb{$Z(u,\cd)$ is $\dbF$-adapted}\\
\ns\ds\qq\qq\qq\qq\qq\q u\in[s,t], \ \|Z(\cd\,,\cd)\|_{L^2_\dbF(s,t;L^2(s,t;H))}^2\equiv
\dbE\int_s^t\int_s^t|Z(u,v)|^2dvdu<\infty\Big\}.\ea$$
We also need two more spaces for operator-valued processes. To this end, we introduce Banach space $\dbB:=\dbR^n\times L^4(0,T;\dbR^n)$.
Obviously $\dbB$ is separable, and there exists numerable dense subset $\dbB_0$ and
$$\dbB_1:= \Big\{\sum_{i=1}^nk_i \a^{i};\ k_i\in\dbQ,\ \a^{i}\in \dbB_0,\ n\in\dbN\Big\}\subset\dbB,$$
where $\dbQ$, $\dbN$ is respectively the set of rational number, integer number in $\dbR$. We denote $\dbB':= \dbR^n\times L^{\frac 4 3}(0,T;\dbR^n)$ the dual space of $\dbB$,  $\cL(H';H'')$ the space of linear bounded operators between Banach spaces $H'$ and $H''$,
\bel{Definition-sL-i}\ba{ll}
\ns\ds \sL_1:=L^2_{1,\cF}\big(0,T;\cL(\dbB;L^{\frac 4 3}(0,T;\dbR^n))\big)\\
\ns\ds \qq=\Big\{\cA:[0,T]\times\Omega\mapsto \cL(\dbB;L^{\frac 4 3}(0,T;\dbR^n))\bigm|\hb{for any}\ b\in\dbB,\ a(\cd)\in L^{4}(0,T;\dbR^n),\\
\ns\ds\qq\q \ \int_0^T\!a(s)^{\top}\!\big[\cA(\cd,\cd)b\big](s)ds\!\in\! L^2_{\dbF}(0,T;\dbR),\ \sup_{t\in[0,T]}\dbE\|\cA(t)\|^2_{\cL(\dbB;L^{\frac 4 3}(0,T;\dbR^n))}\!<\!\infty\!\Big\},\\
\ns\ds \sL_2:=L^2_{1,\cF}\big(0,T;\cL(\dbB;\dbR^n)\big)\\
\ns\ds \qq=\Big\{\cA:[0,T]\times\Omega\mapsto \cL(\dbB;\dbR^n)\bigm|\sup_{t\in[0,T]}\dbE\|\cA(t)\|^2_{\cL(\dbB; \dbR^n)}<\infty,\ \hb{and for any} \ a(\cd)\in\dbR^n, \\
\ns\ds\qq\q \ b\in\dbB,\ a^{\top}\big[\cA(\cd,\cd)b\big]\in L^2_{\dbF}(0,T;\dbR) \Big\}.
\ea
\ee
In this paper, $K$ is a generic constant which varies in different context.

\subsection{Problem formulation}

For FSVIE
(\ref{FSVIE-3-1}), we introduce the following assumptions.

\medskip

{\bf(H1)} Suppose $\varphi(\cd)\in C_{\dbF}([0,T];L^{p}(\Omega;\dbR^n))$, given $p>2$, nonempty $U\subset\dbR^m$, $u(\cd)\in\cU^{p}$, where
$$\cU^{p}:=\Big\{u(\cdot):[0,T]\times\Omega\rightarrow U \ \hb{is measurable and $\dbF$-adapted s.t.} \sup_{t\in[0,T]}\dbE|u(t)|^p<\infty\Big\},$$
$
b,\si:[0,T]^2\times\dbR^n\times U\times\Omega\to\dbR^n$
are measurable such that
$s\mapsto\big(b(t,s,x,u),\si(t,s,x,u)\big)$
is $\dbF$-adapted, $ b,\ \si$ are linear growth of $x,\ u$, twice continuously differentiable of $x$ with bounded first, second order derivatives,
$$\ba{ll}
\ns\ds|f(t,s,x,u)-f(t',s,x,u)| \le\rho(|t-t'|)\big[1+|x|+|u|\big],\
t,t',s\in[0,T],\   f:= b, \si, b_x, \si_x,
\ea$$
with $\rho:[0,\infty)\to[0,\infty)$ a modulus of continuity (continuous, monotone, increasing
function with $\rho(0)=0$). Moreover, $(s,x,u)\!\mapsto \! b_{xx}(t,s,x,u), \si_{xx}(t,s,x,u)$ are continuous uniformly in $t\in[0,T]$.

\medskip

The next result is concerned with the well-posedness of (\ref{FSVIE-3-1}), the proof of which is straightforward adaptation of the counterparts in e.g. \cite{Shi-Wang-Yong 2015}, \cite{Wang-Yong 2013}.

\mds

\bl\label{Lemma-SVIEs}
 Let {\rm(H1)} hold. Then there exists $X(\cd)\in
C_{\dbF}([0,T];L^p(\Omega;\dbR^n))$  satisfying (\ref{FSVIE-3-1})
such that for absolute constant $K$,
\bel{SVIE-estimate-1}\ba{ll}
\ns\ds \sup_{t\in[0,T]}\dbE|X(t)|^p \le K\Big[\sup_{t\in[0,T]}\dbE|\varphi(t)|^p+\sup_{t\in[0,T]}\dbE\int_0^t |b(t,s,0,u(s))|^pds\\
\ns\ds\qq\qq \qq\qq+
\sup_{t\in[0,T]}\dbE\int_0^t |\si(t,s,0,u(s))|^pds\Big].
\ea\ee
Moreover, for $i=1,2$, if $X_i(\cd)$ is the solution associated with $(\varphi_i,b_i,\si_i)$, then
\bel{SVIE-estimate-2}\ba{ll}
\ns\ds \sup_{t\in[0,T]}\dbE|X_1(t)-X_2(t)|^p\\
\ns\ds \le K \Big\{\sup_{t\in[0,T]}\dbE|\varphi_1(t)-\varphi_2(t)|^p +\sup_{t\in[0,T]}\dbE\Big[\int_0^T|b_1(t,s,X_2(s),u(s))-b_2(t,s,X_2(s),u(s))|ds\Big]^p\\
\ns\ds\q+\sup_{t\in[0,T]}\dbE\Big[\int_0^T|\si_1(t,s,X_2(s),u(s))-\si_2(t,s,X_2(s),u(s))|^2ds\Big]^{\frac p 2}\Big\}.
\ea\ee
\el
In the following we define $p:=4+\kappa$ with constant $\kappa>0$, and $\cU_{ad}:= \cU^{4+\kappa}$.

For the involved functions $h$ and $l$ in (\ref{cost3.1}), we make the following assumption.

\bigskip

{\bf (H2)} Let
$
h:\dbR^n\times \Omega\to\dbR$, $l:[0,T]\times\dbR^n\times U\times\Omega\to\dbR$
be measurable such that
$x\mapsto h(x)$, $(x,u)\mapsto l(s,x,u)$
are twice continuously differentiable with
$$\ba{ll}
\big[|h_x(x)|+|l_x(t,x,u)|\big]\leq L[1+|x|+|u|],\ \ \big[|h_{xx}(x)|+|l_{xx}(t,x,u)|\big]\leq L,\ \ x\in\dbR^n,\ \ u\in U.
\ea
$$

We state the optimal control problem as follows.

\mds

\bf Problem (C). \rm Given (\ref{FSVIE-3-1}), (\ref{cost3.1}), we are aiming to
find $\bar u(\cd)\in\cU_{ad}$ such that
$J(\bar u(\cd))=\inf_{u(\cd)\in\cU_{ad}}J(u(\cd)).$

\mds

In above, we call $\bar u$ the $optimal$ $control$, $\bar X$ the $optimal$ $state$ $process$, $(\bar X,\bar u)$ the $optimal$ $pair$.

\section{Maximum Principle for controlled SVIEs}

This section is devoted to obtaining optimality necessary conditions of \bf Problem (C) \rm via spike variation. Without further statement, let $u\in U$, $\t\in[0,T)$, $\e>0$ be sufficiently small such that $\t+\e\leq T$,
\bel{First-definition-E-tau}\ba{ll}
\ns\ds
E_{\t,\e}:=[\t,\t+\e],\ \ u^{\e}(\cd):= uI_{E_{\t,\e}}(\cd)+ \bar u(\cd)I_{[0,T]/E_{\t,\e}}(\cd).
\ea\ee

\mds

\subsection{First-order variational equations and related quadratic functional}

In this part, we derive one quadratic functional of the solutions for first-order variational equations. How to treat this functional appropriately is the crucial step in establishing the maximum principles. These procedures are not necessary with convex control region (\cite{Shi-Wang-Yong 2015}, \cite{Yong 2006}, \cite{Yong 2008}).

Inspired by the SDEs case in e.g. \cite{Peng90}, \cite{Yong-Zhou 1999}, we introduce the first-order variational equations which includes two linear SVIEs as follows,
\bel{variational-equations}\left\{\ba{ll}
\ns\ds X_1(t)= \int_0^t\bar b_x(t,s)X_1(s)ds+\int_0^t\big[\bar\sigma_x(t,s)X_1(s)+\delta\sigma(t,s)\big]dW(s),\\
\ns\ds X_2(t)=\varphi_2(t)+\int_0^t\bar b_x(t,s)X_2(s)ds+\int_0^t \bar\sigma_x(t,s)X_2(s) dW(s).
\ea\right.\ee
Here $t\in[0,T]$, and for $f:= b,\ \si,\ \si_x,$ we make the following conventions,
\bel{Important-notations-1}\left\{\!\!\!\ba{ll}
\ns\ds \varphi_2(t):=\!\!\int_0^t\!\!\big[\frac 1 2\bar b_{xx}(t,s)X_1^2(s)\!+\!\d b(t,s)\big]ds
\!\!+\!\!\int_0^t\!\!\big[\frac 1 2 \bar\si_{xx}(t,s)X_1^2(s)\!+\!\d\si_x(t,s)X_1(s)\big]dW(s),\\
\ns\ds \bar f_x(t,s):= f_x(t,s,\bar X(s),\bar u(s)),\  \d f(t,s):= f(t,s,\bar X(s),u^\e(s))-f(t,s,\bar X(s),\bar u(s)),\\
\ns\ds \bar f_{xx}(t,s)X_1^2(s):= \Big(\hb{tr}\{f_{xx}^1(t,s)X_1(s)X_1(s)^{\top}\},\cdots,\hb{tr}\{f_{xx}^n(t,s)X_1(s)
X_1(s)^{\top}\}\Big)^{\top}.
\ea\right.\ee
We give the following standard estimates, the proof of which is similar to the SDEs case.

\mds

\bl\label{Lem-estimate-Appendix}
Suppose (H1) is true, $X_1, X_2$ satisfy (\ref{variational-equations}), $(\bar X,\bar u)$ is an optimal pair, $X^\e$ is the state process associated with $u^\e$, $u^\e$ is defined in (\ref{First-definition-E-tau}). Then
$$\ba{ll}
\ns\ds \sup_{t\in[0,T]}\dbE |X_1(t)|^2\leq K\e,\ \ \sup_{t\in[0,T]}\dbE|X^\e(t)-\bar X(t)-X_1(t)-X_2(t)|^2\leq o(\e^2).
\ea$$
\el
From Lemma \ref{Lem-estimate-Appendix}, one sees that
\bel{variational-inequality-1}\ba{ll}
\ns\ds o(\e)\leq \dbE\int_0^T\bar l_x(s)^{\top}\big[X_1(s)+X_2(s)\big]ds+\dbE\big[\bar h_x(T)^{\top}(X_1(T)+X_2(T))\big]\\
\ns\ds\qq+\dbE\Big[\int_0^T\d l(s)ds+\frac 1 2 \int_0^TX_1(s)^{\top}\bar l_{xx}(s)X_1(s)ds+\frac 1 2 \big[X_1(T)^{\top}\bar h_{xx}(T)X_1(T)\big]\Big],
\ea\ee
where for example,
\bel{Important-notation-2}\left\{\ba{ll}
\ns\ds \bar l_x(s):= l_x(s,\bar X(s),\bar u(s)),\ \ \bar l_{xx}(s):= l_{xx}(s,\bar X(s),\bar u(s)),\ \
\bar h_x(T):= h_x(\bar X(T)),\\
\ns\ds \bar h_{xx}(T):= h_{xx}(\bar X(T)),\ \  \d l(s):= l(s,\bar X(s),u^\e(s))-l(s,\bar X(s),\bar u(s)),\ \ s\in[0,T].
\ea\right.\ee
We introduce first-order adjoint equation of the form
\bel{first-order-adjoint-equation}\ba{ll}
\ns\ds \bar Y(t)=\bar l_x(t)^{\top}+\bar b_x(T,t)^{\top}\bar h_x(T)+\bar\si_x(T,t)^{\top}\bar\pi(t)+\int_t^T \bar b_x(s,t)^{\top}\bar Y(s)ds\\
\ns\ds\qq\q+\int_t^T\bar\si_x(s,t)^{\top}\bar Z(s,t)ds-\int_t^T \bar Z(t,s)dW(s),
\ea\ee
and the Hamiltonian function
\bel{Hamiltonian-function}\ba{ll}
\ns\ds H(t,x,\bar X(\cd),\bar Y(\cd),\bar Z(\cd,t),u):=
\dbE_t\big[b(T,t,x,u)^{\top}\bar h_x(T)+\si(T,t,x,u)^{\top}\bar\pi(t)\big]\\
\ns\ds\qq\q+l(t,x,u)+\dbE_t\int_t^Tb(s,t,x,u)^{\top}\bar Y(s)ds+\int_t^T\si(s,t,x,u)^{\top}\bar Z(s,t)ds.
\ea\ee
Observe that (\ref{first-order-adjoint-equation}) is a linear backward stochastic Volterra integral equation (BSVIE) which admits a unique pair of solution $(\bar Y(\cd),\bar Z(\cd,\cd))\in L^2_{\dbF}(0,T;\dbR^n)\times  L^2(0,T;L^2_\dbF(0,T;\dbR^{n\times n}))$ such that
$$\ba{ll}
\ns\ds
\bar Y(t)=\dbE_s \bar Y(t)+\int_s^t \bar Z(t,r)dW(r),\ \ s\in[0,t],\ \ t\in[0,T]. \ \ a.e.
\ea$$
Above $(\bar Y,\bar Z)$ is named as M-solutions of BSVIEs, see e.g. \cite{Shi-Wang-Yong 2015}, \cite{Wang-Yong 2013}, \cite{Yong 2008}. Thanks to Lemma \ref{Lem-estimate-Appendix}, we have
\bel{d-si-T-s-Z-X-1}\ba{ll}
\ns\ds \dbE\int_0^T\Big[\d\si_x(T,s)\bar\pi(s)+\int_s^T\d\si_{x}(t,s)\bar Z(t,s)dt\Big]X_1(s)ds=o(\e).
\ea\ee
To sum up, by (\ref{d-si-T-s-Z-X-1}) and Theorem 5.1 in \cite{Yong 2008}, we can transform (\ref{variational-inequality-1}) into,
\bel{quadratic-form-X-1}\ba{ll}
\ns\ds o(1)\!\leq \! \frac 1 \e \dbE\!\int_0^T\! \D H^\e(t)dt+\sE(\e),
\ea\ee
where $ \D H^{\e}(\cd), \ \sE(\e) $ are defined as
\bel{bar-H-xx-definition}\ba{ll}
\ns\ds \D H^\e(t):=H(t,\bar X(t),\bar X(\cd),\bar Y(\cd),\bar Z(\cd,t),u^\e(t))
\!-\!H(t,\bar X(t),\bar X(\cd),\bar Y(\cd),\bar Z(\cd,t),\bar u(t)),\\
\ns\ds \sE(\e):=\frac {1}{2\e} \dbE\Big\{\int_0^T\hb{tr}\big[\bar H_{xx}(t)X_1(t)X_1(t)^{\top}\big]dt+ \hb{tr}\big[\bar h_{xx}(T)X_1(T)X_1(T)^{\top}\big] \Big\},\\
\ns\ds \bar H_{xx}(t):= H_{xx}(t,\bar X(t),\bar X(\cd),\bar Y(\cd),\bar Z(\cd,t),\bar u(t)),\ \ t\in[0,T].
\ea\ee
In addition, from $(\bar Y(\cd),\bar Z(\cd,\cd))$ in (\ref{first-order-adjoint-equation}), one has $\bar H_{xx}(\cd)\in L^2_{\dbF}(0,T;\dbR^{n\times n})$. Therefore, to give the maximum principle, we need to deal with the quadratic form $\sE(\e)$ in (\ref{quadratic-form-X-1}).

To see the encountered challenges, we recall the corresponding procedures of treating $\sE(\e)$ in SDEs case (e.g. \cite{Peng90}, \cite{Yong-Zhou 1999}):
proving some estimates for the solutions of variational equations, deriving the controlled linear SDE satisfied by $X_1(\cd)X_1(\cd)^{\top}$, introducing suitable adjoint equations, and using duality tricks between the forward, backward systems by the well-known It\^{o} formula.
Nevertheless, if we follow above techniques in our setting, we immediately meet some fundamental difficulties which actually indicate the distinctions between the two optimal control problems. The first one lies in the inadequate role of It\^{o} formula in deriving appropriate equation of $X_1(\cd)X_1(\cd)^{\top}$. One may differentiate $X_1$ by imposing differentiability conditions on $b$, $\sigma$ (see Section 2 of \cite{Agram-Oksendal 2015}) and thus makes It\^{o} formula go through. However, this would cause more complicated double integrals with respect to Lebesgue integral and It\^{o} integral.
The second one is concerned with the introducing of proper second-order adjoint equations and suitable duality tricks. According to \cite{Wang-2016-submitted} (see also Remark \ref{Some-essential-ideas}), even in special linear quadratic framework, it seems impossible to construct one $complete$ second-order adjoint equation which directly covers that in SDEs case. In other words, in SVIEs case we need to introduce other essential notions. As to the duality, it surely can not be realized without overcoming the obstacles aforementioned.

To provide more fundamental ideas, we revisit the particular SDEs case from new viewpoints. We consider the second-order adjoint equation, i.e., a linear BSDE of
\bel{Second-order-adjoint-equations-SDEs}\left\{\ba{ll}
\ns\ds dP_2(t)=-\big[\bar b_x(t)^{\top}P_2(t)+P_2(t)\bar b_x(t)+\bar\si_x(t)^{\top}\L_2(t)+\L_2(t)\bar\si_x(t)\\
\ns\ds\qq\qq\qq+\bar H_{xx}(t)+\bar\si_x(t)^{\top}P_2(t)\bar\si_x(t)\big]dt+\L_2(t)dW(t),\ \ t\in[0,T],\\
\ns\ds P_2(T)=\bar h_{xx}(T).
\ea\right.\ee
In some existing literature (\cite{Du-Meng 2013}, \cite{Lv-Zhang 2015}), $P_2$ is called the $second$-$order$ $adjoint$ $process.$ The maximum principle of SDEs (\cite{Peng90}, \cite{Yong-Zhou 1999}) includes $\sH_p:= \d\bar\si(\t)^{\top}P_2(\t)\d\bar\si(\t)$, where
$$\d\bar\si(\t):= \big[\si(\t,\bar X(\t),u)-\si(\t,\bar X(\t),u(\t))\big].$$
Notice that $\sH_p$ is limit counterpart of $\sE(\e)$, as $\e\rightarrow0$.

In the following, we introduce another way to obtain $P_2(\cd)$ without (\ref{Second-order-adjoint-equations-SDEs}). To this end, for any $\t\in[0,T],$ $\xi_i\in L^2_{\cF_\t}(\Omega;\dbR^n)$, $i=1,2$, we define
\bel{xi-i-xi-quadratic-form-SDEs} \left\{\!\!\!\!\ba{ll}
\ns\ds Y_i(t)=\xi_i+\int_\t^t\bar b_x(s)Y_i(s)ds+\int_\t^t\bar \si_x(s)Y_i(s)dW(s),\ \  i=1,2,\\
\ns\ds J(\t,\xi_1,\xi_2):=\dbE_\t\int_\t^T Y_1(s)^{\top}\bar H_{xx}(s)Y_2(s)ds+\dbE_\t\big[Y_1(T)^{\top}\bar h_{xx}(T)Y_2(T)\big].
\ea\right.\ee
By using It\^{o} formula to $Y_1^{\top}P_2Y_2$ on $[\t,T]$ and recalling above $\sH_p$, we see that
$$\xi_1^{\top}P_2(\t)\xi_2=J(\t,\xi_1,\xi_2),\ \ \sH_p=J(\t,\d\bar\si(\t),\d\bar\si(\t)).$$

We take a closer look at $J(\t,\xi_1,\xi_2)$. Notice that the conventional functional analysis theories tell us: given a bounded bilinear operator $F$ on Hilbert space $\cH\times\cH$, there exists a unique operator $G$ on $\cH$ such that $F(x,y)=\lan Gx,y\ran_{\cH}$. Inspired by this point, in the next Lemma \ref{SMP-step-1-SDE-main} we prove that there exists a unique measurable, continuous, $\dbF$-adapted, $\dbR^n$-valued process $\cB_3(\cd)$ such that $J(\t,\xi_1,\xi_2)=\xi_1^{\top}\cB_3(\t)\xi_2$.
The arbitrariness of $\xi_i$ and the continuity of $P_2$, $\cB_3$ lead to %
$$\dbP\big\{\cB_3(\t)=P_2(\t),\ \forall \t\in[0,T]\big\}=1.$$

We observe that the classical maximum conditions only directly relate to $P_2(\cd)$, but not $\L_2(\cd)$. Consequently, above proposed procedures indicate another approach to derive maximum principle without second-order adjoint equation (\ref{Second-order-adjoint-equations-SDEs}). Moreover, one can drop the reliance on It\^{o} formula and the system for $X_1(\cd)X_1(\cd)^{\top}$.
These points provide us the key clues for following-up investigations on SVIEs.

\subsection{Representations of some quadratic functionals}

Given optimal pair $(\bar X(\cd),\bar u(\cd))$, $u\in U$, we define
\bel{d-bar-si-1}
\ba{ll}
\ns\ds
\d\bar\si(t,\t):= \si(t,\t,\bar X(\t),u)-\si(t,\t,\bar X(\t),\bar u(\t)),\\
\ns\ds \D\bar\si(\cd,\t):=\big(\d\bar\si(T,\t),\d\bar\si(\cd,\t)\big),\ \ t,\ \t\in[0,T].
\ea\ee
Under (H1) with $p=4+\kappa$, one has
$\d\bar\si(\cd,\t)\in C_{\dbF}([0,T];L^{4+\kappa}(\Omega;\dbR^n))$.

Given $\bar b_x(\cd)$, $\bar\si_x(\cd)$, $\bar h_{xx}(T)$, $\bar H_{xx}(\cd)$ in (\ref{Important-notations-1}), (\ref{Important-notation-2}), (\ref{bar-H-xx-definition}), similar as (\ref{xi-i-xi-quadratic-form-SDEs}), we introduce
\bel{F-d-bar-si}\left\{\!\!\!\ba{ll}
\ns\ds F^{\d\bar\si,\d\bar\si}(\t):= \dbE_{\t}\int_{\t}^T \dbX (s)^{\top}\bar H_{xx}(s)\dbX (s)ds+\dbE_\t\big[\dbX (T)^{\top}\bar h_{xx}(T)\dbX (T)\big],\\
\ns\ds \dbX (t)\!=\!\d\bar\si(t,\t)\!+\!\int_{\t}^t\bar b_x(t,s)\dbX (s)ds\!+\!\int_{\t}^t\bar\si_x(t,s)\dbX(s)dW(s),\ \ \forall t\in[\t,T].
\ea\right.\ee
Thanks to (H1) and Lemma \ref{Lemma-SVIEs}, there exists a unique $\dbX(\cd)\in C_{\dbF}([\t,T];L^{4+\kappa}(\Omega;\dbR^n))$ satisfying (\ref{F-d-bar-si}).
To represent $F^{\d\bar\si,\d\bar\si}(\t)$, a quadratic functional with respect to $\dbX(\cd)$ for any fixed $\t$, we state the following result.

\mds

\begin{lemma}\label{SMP-step-1-main}
 Suppose (H1), (H2) hold true with $p=4+\kappa$, $(\bar X(\cd),\bar u(\cd))$ is optimal. Then
\bel{main-result-1}\ba{ll}
\ns\ds F^{\d\bar\si,\d\bar\si}(\t)\!=\!\d\bar\si(T,\t)^{\top}\!\big(\cB_1(\t)\D\bar\si(\cd,\t)\big)
\!+\!\!\int_0^T\!\d\bar\si(s,\t)^{\top}\!\big(\cB_2(\t)\D\bar\si(\cd,\t)\big)(s)ds,
\ea\ee
where $\cB_1\in \sL_2$ and $\cB_2\in\sL_1$ satisfy (see (\ref{Definition-sL-i}) for the definitions of $\sL_i$)
\bel{boundedness-db-B-2-1-1-1}\ba{ll}
\ns\ds \Big[\|\cB_1(\t)\|_{\cL(\dbB;\dbR^n)}+\|\cB_2(\t)\|_{\cL(\dbB;L^{\frac 4 3}(0,T;\dbR^n))}\Big]\\
\ns\ds\ \leq K  \Big[\dbE_{\t}\int_\t^T|\bar H_{xx}(s)|^2ds\Big]^{\frac 1 2}+K\Big[\dbE_{\t}|\bar h_{xx}(T)|^2\Big]^{\frac 1 2},\ \ a.s.\ \ \t\in[0,T].
\ea\ee
If there exists $(\wt\cB_1,\wt\cB_2)$ satisfying (\ref{main-result-1}), (\ref{boundedness-db-B-2-1-1-1}), then
\bel{uniqueness-of-cB-1-2}\ba{ll}
\ns\ds
\dbP\big(\o\in\Omega;\ \wt\cB_1(\t,\o)=\cB_1(\t,\o)\big)=\dbP\big(\o\in\Omega;\  \wt\cB_2 (\t,\o)=\cB_2(\t,\o)\big)=1.
\ea\ee
\end{lemma}
\br
In above, we introduce $\cB_1,\ \cB_2$ to treat $\d\bar\si(\cd,\t)$. These two processes are indispensable, and independent with each other (Subsection 3.4.1). If $\d\bar\si(\cd,\t)\equiv \d\bar\si(\t) $, they will be unified into one $\dbR^{n\times n}$-valued process (Lemma \ref{Lemma-special-sigma}).
\er

\br
For almost $\o\in\Omega$, $\t\in[0,T]$, $\d\bar\si(\cd,\t,\omega)\in C([0,T];\dbR^n).$ However, we extend $C([0,T];\dbR^n)$ into $L^4(0,T;\dbR^n)$ since the dual space of the later is easier to treat. This illustrates $L^{\frac{4}{3}}(0,T;\dbR^n)$ in $\cB_2$.
\er

\mds

To prove Lemma \ref{SMP-step-1-main}, for $\a:= (\a_1,\a_2(\cd))\in \dbB$, $\t\in[0,T]$, $t\in[\t,T],$ consider
\bel{auxiliary-state-equation-1}\left\{\ba{ll}
\ns\ds X^{\a}(t)=\a_2(t)+\int_{\t}^tA(t,s)X^{\a}(s)ds+\int_{\t}^tB(t,s)X^{\a}(s)dW(s),\
\ \ \ a.e. \\
\ns\ds \cX^{\a}(T)=\a_1+\int_{\t}^TA(T,s)X^{\a}(s)ds+\int_{\t}^TB(T,s)X^{\a}(s)dW(s).
\ea\right.\ee
Moreover, for $\bar \a,\ \wt \a\in \dbB$, $ \t\in[0,T]$, we define
\bel{auxiliary-cost-1}\ba{ll}
\ns\ds  f_1^{\bar \a,\wt \a}(\t):=\dbE_{\t}\int_{\t}^T X^{\bar \a}(s)^{\top}Q(s)X^{\wt \a}(s)ds+\dbE_\t\big[\cX^{\bar \a}(T)^{\top}G\cX^{\wt \a}(T)\big].\ \ a.s.
\ea\ee

\bigskip

{\bf (H3)} $Q \in L^2_{\dbF}(0,T;\dbR^{n\times n})$, $G\in L^2_{\cF_T}(\Omega;\dbR^{n\times n})$, $A,\ B:[0,T]^2\times\Omega\mapsto\dbR^{n\times n}$ are bounded measurable processes such that for $t\in[0,T]$, $s\mapsto A(t,s),\ B(t,s)$ are $\dbF$-adapted, and with modulus function $\rho(\cd)$,
$$
\big[|A(t,s)|+|B(t,s)|\big]\leq K,\ \ |A(t,s)-A(t',s)|+|B(t,s)-B(t',s)|\leq\rho(|t-t'|),\ \ t,t',s\in[0,T].
$$

\bigskip

\rm
Under (H3), (\ref{auxiliary-state-equation-1}) is solvable with
$$X^{\a}(\cd)\in L^{4}_{\dbF}(\t,T;\dbR^n),\ \ \cX^{\a}(T)\in L^4_{\cF_T}(\Omega;\dbR^n).$$
It is meaningless to discuss $X^{\a}(T)$ if $\a_2(\cd)\in L^4(0,T;\dbR^n)$. Hence we introduce $\cX^{\a}(T)$ and $\a_1\in\dbR^n$ in (\ref{auxiliary-state-equation-1}). The appearance of both $a_2(\cd)$ and $\a_1$ explains the introducing of $\dbB$.

To simplify the notations, we define
\bel{M-Q-G}\ba{ll}
\ns\ds  M^{Q,G}(\t):= \Big[\dbE_{\t}\int_\t^T|Q(s)|^2ds\Big]^{\frac 1 2}+\Big[\dbE_{\t}|G|^2\Big]^{\frac 1 2},\ \ a.s.\ \ \forall\t\in[0,T].
\ea\ee
\bl\label{Lemma-existence-operator-valued}
Suppose (H3) holds. Then for $\a,\ \b\in \dbB$, $\t\in[0,T]$, one has
\bel{auxiliary-1-deterministic}\ba{ll}
\ns\ds f_1^{\a,\b}(\t )=\a_1^{\top}\big(\cB_1(\t )\b\big)+\int_0^T\a_2(t)^{\top}\big(\cB_2(\t )\b\big)(t)dt,\ \  a.s.\
\ea\ee
where $\cB_1\!\in\sL_2$, $\cB_2 \! \in \! \sL_1$, and for $(\t,\o)\in[0,T]\times\Omega$,
\bel{boundedness-db-B}\ba{ll}
\ns\ds \Big[\|\cB_{1}(\t,\o)\|_{\cL(\dbB;\dbR^n)}+ \|\cB_{2}(\t,\o)\|_{\cL(\dbB;L^{\frac 4 3}(0,T;\dbR^n))}\Big]\leq K M^{Q,G}(\t,\o).
\ea\ee
If there is another pair $(\cB_1',\cB_2')$, then for any $\t\in[0,T]$,
\bel{uniqueness-of-cB-1-2-lemma-3.2}\ba{ll}
\ns\ds
\dbP\big(\o\in\Omega;\ \cB_1'(\t,\o)=\cB_1(\t,\o)\big)=\dbP\big(\o\in\Omega;\  \cB_2'(\t,\o)=\cB_2(\t,\o)\big)=1.
\ea\ee
\el

\begin{proof}
\rm  For reader's convenience, we separate the proof into several steps.

\mds

\it Step 1: \rm We  obtain a pair of operator-valued processes on a subset of $[0,T]\times\Omega$ with full measure.

For any $t\in[\t,T]$, $\a\in\dbB$, from (\ref{auxiliary-state-equation-1}) and Gronwall inequality, we see at once that
\bel{Gronwall-SVIE-estimate}\ba{ll}
\ns\ds \dbE_{\t}\Big[\!\int_{\t}^T\!|X^{\a}(s)|^4ds\!+\!|\cX^{\a}(T)|^4\Big]
\!\leq\! K\Big[\int_\t^T|\a_2(s)|^4ds\!+\!|\a_1|^4\Big]\!\equiv \!K \|\a\|_{\dbB}.\
\ea\ee
Consequently, given $\bar\a$, $\wt\a\in\dbB$, $M^{Q,G}(\cd)$ in (\ref{M-Q-G}), it follows that
\bel{bounded-f-1}\ba{ll}
\ns\ds  |f_1^{\bar \a,\wt \a}(\t)|\leq K M^{Q,G}(\t) \|\bar \a\|_{\dbB}  \|\wt \a\|_{\dbB}.\qq a.s.
\ea\ee
For $\bar \a^{i}, \wt \a^{i}, \bar \a, \wt \a\in \dbB_1$, $k, l\in\dbQ$, we define $N\!\equiv \! N(\bar \a^{i},\ \wt \a^{i},\bar \a,\wt \a,k,l,Q,G)$ and $\cN$ as,
\bel{auxiliary-N-1}\left\{\!\!\!\!\{\ba{ll}
\ns\ds N:= \Big\{\!(\t,\o)\in[0,T]\times\Omega; |f_1^{\bar \a,\wt \a}(\t,\o)|\!\leq \! KM^{Q,G}(\t,\o)\|\bar \a\|_{\dbB}\cd\|\wt \a\|_{\dbB}, \\
\ns\ds\qq \qq f_1^{k\bar \a^{1}+l\bar \a^{2},\wt \a}(\t,\o)=k f_1^{\bar \a^1,\wt \a}(\t,\o)+ l f_1^{\bar \a^2,\wt \a}(\t,\o),\\
\ns\ds\qq\qq
f_1^{\bar \a, k\wt \a^{1}+l\wt \a^{2}}(\t,\o)=k f_1^{\bar \a,\wt \a^1}(\t,\o)+ l f_1^{\bar \a,\wt \a^2}(\t,\o)\Big\},\\
\ns\ds \cN\equiv\cN(G,Q):= \bigcap_{\bar \a,\wt\a\in \dbB_1}\bigcap_{\bar \a^i,\wt \a^i\in \dbB_1}\bigcap_{k,l\in\dbQ}N(\bar \a^{i},\ \wt \a^{i},\bar \a,\wt \a,k,l,Q,G).
\ea\right.\ee
Notice that $\big[\l\times\dbP\big](N)=T$, $\big[\l\times\dbP\big](\cN)=T$, where $\l$ is the Lebesgue measure. In addition, by inequality (\ref{bounded-f-1}), for any $\t\in[0,T]$, one has $\dbP(\cN_{\t})=1$ with $\cN_{\t}:= \big\{\omega\in\Omega;\ (\t,\omega)\in\cN\big\}$.

For any $(\t,\o)\in\cN$, it is easy to see that $f_1^{\cd,\cd}(\t,\o):\dbB_1\times \dbB_1\mapsto\dbR$ is a bounded bilinear map in the sense of (\ref{auxiliary-N-1}). According to Lemma \ref{Lemma-existence-operator-process-Appendix}, there exists a unique linear bounded functional $\cB_{1,1}(\t,\o): \dbB\mapsto\dbR^n$ and a unique linear bounded operator $\cB_{1,2}(\t,\o):\dbB\mapsto L^{\frac 4 3}(0,T;\dbR^n)$ such that for $\a,\ \b\in \dbB_1$,
\bel{cB-1-1-cB-1-2}\ba{ll}
\ns\ds f_1^{\a,\b}(\t,\o)=\a_1^{\top}\big(\cB_{1,1}(\t,\o)\b\big)+\int_0^T\a_2(t)^{\top}\big(\cB_{1,2}(\t,\o)\b\big)(t)dt.\ \
\ea\ee
For any $\t\in[0,T]$, recalling $\dbP(\cN_\t)=1$, we know that  (\ref{cB-1-1-cB-1-2}) holds almost surely.

\mds

\it Step 2: \rm We deduce a pair of operator-valued processes on $[0,T]\times\Omega$.

Given $\a,\b\in\dbB$, there exist $\{\a_n\}_{n=1}^{\infty},\ \{\b_n\}_{n=1}^{\infty}\subset\dbB_1$ such that
$$\big[\|\a_n-\a\|^4_{\dbB}+\|\b_n-\b\|^4_{\dbB}\big]\rightarrow0,\ \ n\rightarrow\infty.$$
For any $(\t,\o)\in\cN$, we denote by $\sB_n(\t,\o)$, $\sB(\t,\o)$
$$\ba{ll}
\ns\ds
\sB_n(\t,\o):= \a_{n,1}^{\top}\big(\cB_{1,1}(\t,\o)\b_n\big)
+\int_0^T\a_{n,2}(t)^{\top}\big(\cB_{1,2}(\t,\o)\b_n\big)(t)dt,\\
\ns\ds \sB(\t,\o):= \a_1^{\top}[\cB_{1,1}(\t,\o)\b]+\int_0^T\a_2(t)^{\top}\big(\cB_{1,2}(\t,\o)\b\big)(t)dt.
\ea
$$
One has $\lim\limits_{n\rightarrow\infty}\Big|\sB(\t,\o)-
\sB_n(\t,\omega)\Big|=0.$
On the other hand, by (\ref{cB-1-1-cB-1-2}), $\sB_n(\cd)$ is measurable, adapted process. Hence similar conclusion holds for $\sB(\cd)$.

At this moment, we define two processes $ \cB_i(\t,\o)$ on $[0,T]\times\Omega$ as
\bel{Definition-of-cB-i}\ba{ll}
\ns\ds \cB_1(\t,\o)\!:=\!\cB_{1,1}(\t,\o)I_{\cN}(\t,\o),\
\cB_2(\t,\o)\!:=\!\cB_{1,2}(\t,\o)I_{\cN}(\t,\o).
\ea\ee
For any $\a_1\in\dbR^n$, by choosing $\a:=(\a_1,0)$ in $\sB(\cd)$, it follows from (\ref{Definition-of-cB-i}) that
$$\ba{ll}
\ns\ds
\a_1^{\top}\big[\cB_1(\t,\o)\b\big]\!=\!\sB(\t,\o)I_{ \cN}(\t,\o)\! =\!\Big[\!\lim\limits_{n\rightarrow\infty}\sB_{n}(\t,\o)\!\Big]I_{ \cN}(\t,\o).
\ea$$
It is then evident to deduce the measurability of $(\t,\o)\mapsto\a_1^{\top}\big(\cB_1(\t,\o)\b\big)$,
as well as the adaptness.

Similarly one can obtain the case of $\cB_2(\cd,\cd)$.

\mds

\it Step 3: \rm For any $\t\in[0,T]$, $\a,\ \b\in\dbB,$ we prove (\ref{auxiliary-1-deterministic}).

From Step 1, for any $\t\in[0,T]$ and $\omega\in\cN_\t$, (\ref{cB-1-1-cB-1-2}) holds true with any $\a_n,\ \b_n\in\dbB_1$. Similar as (\ref{Gronwall-SVIE-estimate}), for any $\t\in[0,T],$ the following is true almost surely,
$$\ba{ll}
\ns\ds \dbE_\t\int_\t^T|X^{f_n}(t)-X^{f}(t)|^4dt+\dbE_\t|\cX^{f_n}(T)-\cX^{f}(T)|^4\leq K\|f_n-f\|_{\dbB},\ \ f:= \a,\ \b.
\ea$$
As a result, $\lim\limits_{n\rightarrow\infty}\big|f_1^{\a,\b}(\t)-f_1^{\a_n,\b_n}(\t)\big|=0$.
Consequently,
$$\ba{ll}
\ns\ds f_1^{\a,\b}(\t)=\sB(\t ):=\a_1^{\top}\big(\cB_{1,1}(\t )\b\big)
+\int_0^T\a_2(t)^{\top}\big(\cB_{1,2}(\t )\b\big)(t)dt.\ \ a.s.
\ea$$
Our conclusion then follows from the relation between $\cB_{1,i}(\cd,\cd)$ and $\cB_i(\cd,\cd)$ in (\ref{Definition-of-cB-i}).

\mds

\it Step 4: \rm In this part, we discuss the integrability and uniqueness of $\cB_i$.

For $(t,\o)\in\cN$, by Lemma \ref{Lemma-existence-operator-process-Appendix}, we have
$$\ba{ll}
\ns\ds \Big[\|\cB_{1,1}(\t,\o)\|_{\cL(\dbB;\dbR^n)}+ \|\cB_{1,2}(\t,\o)\|_{\cL(\dbB;L^{\frac 4 3}(0,T;\dbR^n))}\Big]\leq K M^{Q,G}(\t,\o).
\ea$$
Hence (\ref{boundedness-db-B}) is lead by (\ref{Definition-of-cB-i}). On the other hand,
from (\ref{bounded-f-1}), $\sup\limits_{\t\in[0,T]}\dbE|f^{ \a,\b}_1(\t)|^2<\infty$.
As a result, for $\a:=(\a_1,0)\in\dbB$, $\a_1\in\dbR$, $\b\in\dbB$, it follows from (\ref{auxiliary-1-deterministic}) that
$\sup\limits_{\t\in[0,T]}\dbE\big|\a_1^{\top}\big(\cB_1(\t)\b\big)\big|^2<\infty.
$

Similarly, we derive the case of $\cB_2(\cd)$.

Eventually, we emphasize that the uniqueness of $\cB_i(\cd)$ is obvious to obtain.
\end{proof}

\mds

The next lemma yields information about an extension of (\ref{auxiliary-1-deterministic}). To this end, for $\t\in[0,T]$ we denote $L^4_{\cF_\t}(\Omega;\dbB)$ the set of $\cF_\t$-strongly measurable $\dbB$-valued random variable $\xi$ such that $\dbE\|\xi\|_{\dbB}^4<\infty$. Recall that a $\dbB$-valued random variable $\xi$ is named $\cF_\t$-strongly measurable if there exists a sequence of $\dbB$-valued simple random variables $\xi_k$ converging to $\xi$.

\mds


\bl\label{Lemma-existence-operator-valued-extension}
For any $\t\in[0,T],$ $\xi:=(\xi_1,\xi_2(\cd)),\ \eta:=(\eta_1,\eta_2(\cd))\in L^4_{\cF_{\t}}(\Omega;\dbB)$,
$$\ba{ll}
\ns\ds f^{\xi,\eta}_1(\t)=\xi_1^{\top}\big[\cB_1(\t)\eta\big]+\int_0^T\xi_2(s)^{\top}\big[\cB_2(\t)\eta\big](s)ds.\qq a.s.
\ea$$
\el

\begin{proof}
\rm
We begin with the simple random variable case. By defining
$$\xi(\o):=\sum_{i=1}^n x_iI_{A_i}(\o),\ \ \eta(\o):=\sum_{j=1}^m y_jI_{B_j}(\o),\ \ \o\in\Omega,\ \ A_i, \ B_j\in\cF_\t,\ \ x_i,\ y_j\in \dbB,
$$
it is easy to see that
$$\ba{ll}
\ns\ds \sum_{i=1}^n \sum_{j=1}^m\Big[x_{i,1}^{\top}[\cB_1(\t)y_j]+\int_0^Tx_{i,2}(s)^{\top}
[\cB_2(\t)y_j](s)ds\Big]\cd I_{A_i}  I_{B_j} \\
\ns\ds =\xi_1^{\top}[\cB_1(\t)\eta]+\int_0^T\xi_{2}(s)^{\top}[\cB_2(\t)\eta](s)ds.
\ea$$
From Lemma \ref{Lemma-existence-operator-valued} we conclude that, for any $ x_i:=(x_{i,1},x_{i,2}(\cd))\in \dbB,\ y_j\in \dbB$,
\bel{f-1-simple-case-0}\ba{ll}
\ns\ds f_1^{x_i,y_j}(\t)=x_{i,1}^{\top}\big(\cB_1(\t )y_j\big)+\int_0^Tx_{i,2}(t)
^{\top}\big(\cB_2(\t )y_j\big)(t)dt.\ \  a.s. \ \
\ea\ee
We thus obtain the desirable conclusion by
$$\ba{ll}
\ns\ds f_1^{\xi,\eta}(\t)
= \sum_{i=1}^n \sum_{j=1}^m\Big(\dbE_\t\int_\t^TX^{x_i}(s)^{\top} Q(s)X^{x_j}(s)ds+
\dbE_\t\big[ \cX^{x_i}(T)^{\top} G \cX^{x_j}(T) \big]\Big)\cd I_{A_i}  I_{B_j} \\
\ea$$
The task now is to treat the general case. For any $\xi:=(\xi_1,\xi_2(\cd)),\ \eta:=(\eta_1,\eta_2(\cd))\in L^4_{\cF_\t}(\Omega;\dbB)$, there exist $\{\xi_n\}_{n\geq1}, \ \{\eta_n\}_{n\geq1}$ such that
$$\|k-k_n\|^4_{L^2_{\cF_\t}(\Omega;\dbB)}\rightarrow0,\ \ n\rightarrow\infty, \ \ k:= \xi,\ \eta.$$
Therefore, similar as (\ref{Gronwall-SVIE-estimate}), when $n\rightarrow\infty$, we have
$$\ba{ll}
\ns\ds \dbE\int_\t^T|X^{k_n}(t)-X^{k}(t)|^4dt+\dbE|\cX^{k_n}(T)-\cX^{k}(T)|^4\leq K\|k-k_n\|^4_{L^2_{\cF_\t}(\Omega;\dbB)}\rightarrow0,
\ea$$
with $k:= \xi,\ \eta.$ This implies that for any $\t\in[0,T]$, $\lim\limits_{n\rightarrow\infty}\dbE|f_1^{\xi,\eta}(\t)-f_1^{\xi_n,\eta_n}(\t)|=0$.

On the other hand, by the estimates in Lemma \ref{Lemma-existence-operator-valued},
$$\ba{ll}
\ns\ds\lim_{n\rightarrow\infty} \dbE\Big|\xi_1^{\top}[\cB_1(\t)\eta]\!-\!\xi_{n,1}^{\top}[\cB_1(\t)\eta_n]\!+\!
\int_0^T\big[\xi_2(t)^{\top}\big(\cB_2(\t)\eta\big)(t)\!-\!
\xi_{n,2}(t)^{\top}\big(\cB_2(\t)\eta_n\big)(t)\big]dt\Big|=0.
\ea$$
Since $(\xi_n,\eta_n)$ are all simple random variables, for any $n\geq1$ we have
$$\ba{ll}
\ns\ds \dbE\Big|f^{\xi_n,\eta_n}_1(\t)- \xi_{n,1}^{\top}\big[\cB_1(\t)\eta_{n}\big]-\int_0^T
 \xi_{n,2}(s)^{\top}\big[\cB_2(\t)\eta_{n}\big](s)ds\Big|=0.
\ea$$
To sum up, the conclusion is followed by
$$\ba{ll}
\ns\ds \dbE\Big|f_1^{\xi,\eta}(\t)-\xi_1^{\top}[\cB_1(\t)\eta]-\int_0^T\xi_2(t)
^{\top}\big(\cB_2(\t)\eta\big)(t)dt\Big|=0.
\ea$$
\end{proof}

Now we give the proof of Lemma \ref{SMP-step-1-main}.

\begin{proof}
\rm Given $\d\bar\si(\cd,\t),\ \D\bar\si(\cd,\t)$ in (\ref{d-bar-si-1}), $\dbX(\cd)$ in (\ref{F-d-bar-si}), like (\ref{auxiliary-state-equation-1}) we introduce
$$\ba{ll}
\ns\ds \cX^{\bar\si}(T):=\d\bar\si(T,\t)+\int_{\t}^T\bar b_x(T,s)\dbX (s)ds+\int_{\t}^T\bar\si_x(T,s)\dbX(s)dW(s).
\ea$$
We also introduce $f_1^{\D\bar\si,\D\bar\si}(\t)$
$$\ba{ll}
\ns\ds
f_1^{\D\bar\si,\D\bar\si}(\t):= \dbE_{\t}\int_{\t}^T \dbX (s)^{\top}\bar H_{xx}(s)\dbX(s)ds+\dbE_\t\big[\cX^{\bar\si}(T)^{\top}\bar h_{xx}(T)\cX^{\bar\si}(T)\big].
\ea$$
By virtue of Lemma \ref{Lemma-existence-operator-valued} and Lemma \ref{Lemma-existence-operator-valued-extension}, there exist $\cB_1(\cd,\cd)$, $\cB_2(\cd,\cd)$ satisfying (\ref{boundedness-db-B-2-1-1-1}) and
$$\ba{ll}
\ns\ds f_1^{\d\bar\si,\d\bar\si}(\t)=\d\bar\si(T,\t)^{\top}\big(\cB_1(\t)\D\bar\si(\cd,\t)\big)
+\int_0^T\d\bar\si(s,\t)^{\top}\big(\cB_2(\t)\D\bar\si(\cd,\t)\big)(s)ds.
\ea$$
The integrability of $\dbX (\cd)$ and $\d\bar\si(\cd,\t)$ yield $\dbE_\t\big|\cX^{\bar\si}(T)-\dbX (T)\big|^4=0$, a.s., $\t\in[0,T]$.

As a result, the conclusion (\ref{main-result-1}) follows immediately.
\end{proof}

\mds

In the rest of this subsection, we discuss the case when $\d\bar\si(\cd,\t)\equiv\d\bar\si(\t) $.

To this end, for any $a_1,\ a_2\in\dbR^n$, $\t\in[0,T]$, similar as (\ref{auxiliary-cost-1}) we define $f^{a_1,a_2}_2(\t)$
\bel{auxiliary-cost-dbR-1}\ \ \ \ba{ll}
\ns\ds  f^{a_1,a_2}_2(\t):=\dbE_\t\int_\t^TX^{a_1}(s)^{\top}Q(s)X^{a_2}(s)ds+\dbE_\t\big[X^{a_1}(T)^{\top} GX^{a_2}(T)\big],\ \ a.s.
\ea\ee
associated with $Q(\cd)$, $G$, where
\bel{auxiliary-state-equation-dbR-1}\qq \ba{ll}
\ns\ds X^{a_i}(t)=a_i+\int_\t^tA(t,s)X^{a_i}(s)ds+\int_\t^tB(t,s)X^{a_i}(s)dW(s),\qq \forall t\in[\t,T].
\ea\ee
In particular, we have $f^{e_i,e_j}_2(\cd)$ where $e_i\in\dbR^n$.

\mds

\bl\label{Lemma-special-sigma}
Suppose (H3) holds true. Then for any $\t\in[0,T],$ $\xi_i\in L^4_{\cF_\t}(\Omega;\dbR^n),$
\bel{boundedness-db-B-3-remark}\ \ \ \ba{ll}
\ns\ds f^{\xi_1,\xi_2}_2(\t)=\xi_1^{\top}\cB_3(\t)\xi_2,\  |\cB_3(\t)|\leq K \Big\{\Big[\dbE_{\t}\int_\t^T|Q(s)|^2ds\Big]^{\frac 1 2}+\Big[\dbE_{\t}|G|^2\Big]^{\frac 1 2}\Big\},
\ea\ee
where $\cB_3(\cd):= \Big\{f^{e_i, e_j}_2(\cd)\Big\}_{1\leq i,j\leq n}$.
If there is another continuous process $ \cB_3'(\cd)$ satisfying (\ref{boundedness-db-B-3-remark}), then
\bel{uniqueness-of-cB-1-3-remark}\ba{ll}
\ns\ds
\dbP\big(\o\in\Omega;\ \cB_3'(\t,\o)=\cB_3(\t,\o),\ \ \forall \t\in[0,T]\big)=1.\ \
\ea\ee
\el

\begin{proof}
 \rm The ideas of this proof are essentially the same as Lemma \ref{SMP-step-1-main} and Lemma \ref{Lemma-existence-operator-valued-extension}. For reader's convenience, we give a sketch as follows.

Given $\bar\a$, $\wt\a\in\dbQ^n$ the set of $n$-dimensional vectors with each component being rational number, and $M^{Q,G}(\cd)$ in (\ref{M-Q-G}), we can deduce that
\bel{bounded-f-2}\ba{ll}
\ns\ds  |f_2^{\bar \a,\wt \a}(\t)|\leq K M^{Q,G}(\t)\cd\|\bar \a\|_{\dbQ^n}\cd \|\wt \a\|_{\dbQ^n}.\qq a.s.
\ea\ee
Moreover, for any $\bar \a^{i},\ \wt \a^{i},\ \bar \a,\ \wt \a\in \dbQ^n$ with $i=1,2,$ and $k,\ l\in\dbQ$, we define
\bel{auxiliary-N-1-real-value}\left\{\!\!\ba{ll}
\ns\ds \wt N:= \Big\{(\t,\o)\in[0,T]\times\Omega;\ |f_2^{\bar \a,\wt \a}(\t,\o)|\leq KM^{Q,G}(\t,\o)\|\bar \a\|_{\dbR^n}\cd\|\wt \a\|_{\dbR^n}, \\
\ns\ds\qq\qq f_2^{k\bar \a^{1}+l\bar \a^{2},\wt \a}(\t,\o)=k f_2^{\bar \a^1,\wt \a}(\t,\o)+ l f_2^{\bar \a^2,\wt \a}(\t,\o),\\
\ns\ds\qq\qq
f_2^{\bar \a, k\wt \a^{1}+l\wt \a^{2}}(\t,\o)=k f_1^{\bar \a,\wt \a^1}(\t,\o)+ l f_2^{\bar \a,\wt \a^2}(\t,\o)\Big\},\\
\ns\ds  \wt\cN := \bigcap_{\bar \a,\wt\a\in \dbQ^n}\bigcap_{\bar \a^i,\wt \a^i\in \dbQ^n}\bigcap_{k,l\in\dbQ}\wt N(\bar \a^{i},\ \wt \a^{i},\bar \a,\wt \a,k,l,Q,G).
\ea\right.\ee
It is easy to see that $\big[\l\times\dbP\big](\wt\cN)=T$ such that for any $(\t,\omega)\in\wt\cN$, $f_2^{\cd,\cd}(\t,\o):\dbQ^n\times \dbQ^n\mapsto\dbR$ is a bounded bilinear map. Moreover, for any $\t\in[0,T]$, $\dbP(\wt\cN_{\t})=1$, where $\wt\cN_\t:= \{\o\in\Omega,(\t,\o)\in\wt\cN\}$.

Using Lemma \ref{Lemma-existence-real-valued-process-Appendix} there exists a unique $\dbR^{n\times n}$-valued matrix $\cB_{1,3}(\t,\o)$, $(\t,\o)\in\wt\cN,$ such that the following holds true with $\a,\ \b\in \dbQ^n$,
\bel{cB-1-1-cB-1-2-2-real}\ba{ll}
\ns\ds f_2^{\a,\b}(\t,\o)=\a^{\top}\cB_{1,3}(\t,\o)\b,\ \ |\cB_{1,3}(\t,\o)|\leq  K M^{Q,G}(\t,\o).
\ea\ee
Let $\a:= e_i$, $\b:= e_j$, $i,\ j=1,\cdots,n$, we obtain that $\cB_{1,3}^{i,j}(\t,\o)=f^{e_i, e_j}_2(\t,\o)$ which is component in the $i$-th line, $j$-th column. In other words, given (\ref{boundedness-db-B-3-remark}), $\cB_{1,3}=\cB_3$ in $\wt\cN$. Moreover,
$$\dbP\big(\o\in\Omega;\cB_{1,3}(\t,\o)=\cB_3(\t,\o)\big)=\dbP(\wt\cN_\t)=1,\ \ \t\in[0,T].$$
Considering (\ref{cB-1-1-cB-1-2-2-real}), we have $\sup\limits_{\t\in[0,T]}\dbE|\cB_3(\t)|^2 <\infty$, and the second result in (\ref{boundedness-db-B-3-remark}).
The measurability, adaptness and continuity of $\cB_3$ are easy to see.

Following the same ideas as in Step 3 of Lemma \ref{Lemma-existence-operator-valued}, and Lemma \ref{Lemma-existence-operator-valued-extension}, for any $\xi_i\in L^4_{\cF_\t}(\Omega;\dbR^n)$, one has $f^{\xi_1,\xi_2}_2(\t)=\xi^{\top}\cB_3(\t)\xi_2$. a.s.

Eventually, (\ref{uniqueness-of-cB-1-3-remark}) follows from the arbitrariness of $\xi\in L^4_{\cF_\t}(\Omega;\dbR^n)$ and continuity of $\cB_3(\cd)$.
\end{proof}

\mds

\br
 Unlike Lemma \ref{Lemma-existence-operator-valued}, several more clearer pictures are given here. The first one is about the better regularity of $ X^{\a}(\cd)$, which saves us from introducing new terms such as $\cX^\a(T)$ of (\ref{auxiliary-state-equation-1}). The second one is about $\dbR^n$-valued process $\cB_3(\cd)$, which plays the same role as operator-valued processes $\cB_1(\cd)$, $\cB_2(\cd)$. Moreover, $\cB_3(\cd)$ has more stronger properties such as measurability, continuity, uniqueness.
\er

\mds

Using Lemma \ref{Lemma-special-sigma}, we give the following result that is comparable with Lemma \ref{SMP-step-1-main},

\mds

\begin{lemma}\label{SMP-step-1-SDE-main}
Suppose (H1), (H2) hold true with $p=4+\kappa$, $\d\bar\si(\cd,\t)\equiv \d\bar\si(\t)$, $(\bar X(\cd),\bar u(\cd))$ is optimal. Then there exist a unique (in the sense of (\ref{uniqueness-of-cB-1-3-remark})) measurable, adapted, continuous $\dbR^{n\times n}$-valued process $\cB_3(\cd):= \big\{F^{e_i, e_j}(\cd)\big\}_{1\leq i,j\leq n}$ such that,
\bel{main-result-1-SDE-special}\ba{ll}
\ns\ds F^{\d\bar\si,\d\bar\si}(\t)=\d\bar\si( \t)^{\top} \cB_3(\t)\d\bar\si( \t),\ \ a.s.\ \ \t\in[0,T].
\ea\ee
\end{lemma}

\subsection{Some subtle asymptotic analyses}

In order to obtain maximum principle, we are in a position to explore some essential relations between $\sE(\e)$ in (\ref{bar-H-xx-definition}) and $F^{\d\bar\si,\d\bar\si}(\cd)$ in (\ref{F-d-bar-si}). To get more intuitive feelings, we look at the case of $n=1$, $\bar b_x(\cd,\cd)=\bar\si_x(\cd,\cd)=0$.
In this case, given deterministic $\cQ(\cd)$, using Fubini theorem and
$$X_1(\cd)=\int_\t^{\cd}\d\bar\si(\cd,s)I_{[\t,\t+\e]}(s)dW(s),$$
we see that
$$\ba{ll}
\ns\ds \frac 1\e \dbE\int_\t^T\cQ(t)\Big|X_1(t)\Big|^2dt
=\frac 1 \e \dbE\int_\t^{\t+\e}\int_s^T\cQ(t)|\d\bar\si(t,s)|^2dtds.
\ea$$
For $\t\in[0,T)$, a.e., by Lebesgue differentiation theorem,
\bel{step-2-simple-case}\ \ \ \ \ \ba{ll}
\ns\ds \lim_{\e\rightarrow0}\frac 1\e \dbE\int_\t^T\cQ(t)\Big|X_1(t)\Big|^2dt=\dbE \int_\t^T\cQ(t)|\d\bar\si(t,\t)|^2dt=\dbE\int_\t^T\cQ(t)|\dbX(t)|^2dt.
\ea\ee
Above (\ref{step-2-simple-case}) indicates certain asymptotic connection between $X_1(\cd)$ and $\dbX(\cd)$. We use this basic idea in our framework and present the following result.

\mds

\bl\label{Theorem-step-2}
  Suppose (H1-H2) hold true, $(\bar X(\cd),\bar u(\cd))$ is optimal pair. Then there exists $\{\e_n\}_{n\geq1}$ such that $\e_n\rightarrow0$ as $n\rightarrow\infty$, and
$$\ba{ll}
\ns\ds \lim_{\e_n\rightarrow0}\frac {1} {\e_n}\Big[\dbE\int_0^T X_1(s)^{\top}\bar H_{xx}(s)X_1(s)ds+\dbE\big[X_1(T)^{\top}\bar h_{xx}(T)X_1(T)\big]\Big]\\
\ns\ds=\dbE\int_\t^T\dbX(s)^{\top}\bar H_{xx}(s)\dbX(s)ds+\dbE\big[\dbX(T)^{\top}\bar h_{xx}(T)\dbX(T)\big].
\ea$$
\el
To prove Lemma \ref{Theorem-step-2}, we take a closer look at $X_1(\cd)$ in (\ref{variational-equations}). Actually, according to its definition, one can rewrite $X_1(\cd)$ as,
\bel{variational-equation-X-1}\ba{ll}
\ns\ds X_1(t)=\left\{\ba{ll}
\ns\ds 0,\qq  t\in[0,\t],\\
\ns\ds \int_\t^t\bar b_x(t,s)X_1(s)ds+\int_\t^t\big[\bar\si_x(t,s)X_1(s)+\d\bar\si(t,s)\big]dW(s),\ t\in[\t,\t+\e],\\
\ns\ds \rho_1(t)+\int_{\t+\e}^t \bar b_x(t,s)X_1(s)ds+\int_\t^t\bar\si_x(t,s)X_1(s) dW(s),\ t\in[\t+\e,\infty),
\ea\right.
\ea\ee
where $\rho _1(\cd)\in C([\t+\e,T];L^4(\Omega;\dbR))$ is $\cF_{\t+\e}$-measurable defined as,
$$\ba{ll}
\ns\ds \rho_1(t):= \int_\t^{\t+\e}\bar b_x(t,s)X_1(s)ds+\int_\t^{\t+\e}\big[\bar\si_x(t,s)X_1(s)+\d\bar\si(t,s)\big]dW(s),\ \ t\geq \t+\e.
\ea$$
To deal with the quadratic form of $X_1(\cd)$ (see Lemma \ref{Lemma-step-2-1}), for small $\e>0$, we introduce $Y_1(\cd)$ on $[\t+\e,T]$,
\bel{Y-1-e}\ba{ll}
\ns\ds Y_1(t)=\varrho_1(t)+\int_{\t+\e}^t\bar b_x(t,s)Y_1(s)ds+\int_{\t+\e}^t\bar\si_x(t,s)Y_1(s)dW(s),
\ea\ee
where
$$\varrho_1(\cd):= \e^{-\frac 1 2}\d\bar\si(\cd,\t)\Big(W(\t+\e)-W(\t)\Big).$$
We claim that $\varrho_1(\cd)\in C_{\cF_{\t+\e}}([\t+\e,T],L^4(\Omega;\dbR^n))$. Therefore, thanks to Lemma \ref{Lemma-SVIEs} and (H1), (\ref{Y-1-e}) admits a unique solution $Y_1(\cd)\in  C_{\dbF}([\t+\e,T],L^4(\Omega;\dbR^n))$.

In fact, since $\bar u(\cd)\in\cU_{ad}$, from Lemma \ref{Lemma-SVIEs} and (H1), we have $\sup\limits_{\t\in[0,T]}\dbE|\bar X(\t)|^{4+\kappa}<\infty$, and
$$\sup\limits_{t\in[0,T]}\dbE|\d\bar\si(t,\t)|^{4+\kappa}\leq L\Big[1+\sup\limits_{\t\in[0,T]}\dbE|\bar X(\t)|^{4+\kappa}+
\sup\limits_{\t\in[0,T]}\dbE|\bar u(\t)|^{4+\kappa}\Big]<\infty.$$
Consequently, by virtue of H\"{o}lder inequality and Jensen's inequality of expectation,
\bel{Brownian-motion-formula-application}\ba{ll}
\ns\ds \!\sup_{t\in[\t+\e,T]}\!\dbE|\varrho_1(t)|^4\!\leq \! \sup_{t\in[\t+\e,T]}\!\Big[\dbE|\d\bar\si(t,\t)|^{4+\kappa}\Big]^{\frac{4}{4+\kappa}}
\e^{-2}\Big|\dbE\big[
|W(\t+\e)\!-\!W(\t)|^{\frac{4(4+\kappa)}{\kappa}}\big]\Big|^{\frac{\kappa}{4+\kappa}}\\
\ns\ds\qq\qq\qq\qq\leq \sup_{t\in[\t+\e,T]}\Big[\dbE|\d\bar\si(t,\t)|^{4+\kappa}\Big]^{\frac{4}{4+\kappa}}\cd \Big[\frac{\big(4[p]+4\big)!}{2^{2[p]+2}\cd\big(2[p]+2\big)!}\Big]^{\frac{p}{[p]+1}}<\infty.
\ea\ee
Here $p:= \frac{4+\kappa}{\kappa}$, $[p]$ is the integer part of $p$, and we used one formula on Brownian motions:
$$\dbE|W(t)|^{2k}=\frac{(2k)!}{2^k k!} t^k,\ \ k\in\dbN,\ \ t\in\dbR^+.$$
Moreover,
$$
\lim\limits_{t\rightarrow t_0}\dbE|\varrho_1(t)-\varrho_1(t_0)|^4=0,\ \ t_0\in[\t+\e,T].
$$
Hence the conclusion of $\varrho_1(\cd)$ is obvious. For natational simplicity, we denote
\bel{Notations-1-dbH}\left\{\!\!\ba{ll}
\ns\ds \dbH(\e,X_1):=\frac 1 {\e}\dbE\Big[\int_{\t+\e}^T X_1(s)^{\top}\bar H_{xx}(s)X_1(s)ds+ X_1(T)^{\top} \bar h_{xx}(T)X_1(T) \Big],\\
\ns\ds \dbH(\e,Y_1):=\dbE\Big[\int_{\t+\e}^T Y_1(s)^{\top}\bar H_{xx}(s)Y_1(s)ds+Y_1(T)^{\top} \bar h_{xx}(T)Y_1(T)\Big].
\ea\right.\ee
Recall that $L^4(0,T;L^4(\Omega;\sX))$ is the set of $\cB([0,T])\otimes\cF_T$-strongly measurable $\sX$-valued process $f(\cd)$ satisfying $\dbE\int_0^T\|f(s)\|_{\sX}^4ds<\infty.$

Using similar ideas as in Lemma 2.5 of \cite{Lv-Zhang 2015}, we present the following result.
\mds

\bl\label{Lemma-Banach-space-result}
Given Banach space $\sX$, suppose $f(\cd)\in L^4(0,T;L^4(\Omega;\sX))$. Then there exists a sequence $\{\e_n\}$ such that,
$$\ba{ll}
\ns\ds \lim_{n\rightarrow\infty}\frac {1}{\e_n}\int_t^{t+\e_n}\dbE\|f(s)-f(t)\|_{\sX}^4ds=0,\qq t\in[0,T].\ \ a.e.
\ea$$
\el
\bl\label{Lemma-step-2-1}
 For $X_1(\cd)$, $Y_1(\cd)$, $\dbH(\e,X_1),\ \dbH(\e,Y_1)$ in (\ref{quadratic-form-X-1}), (\ref{Y-1-e}), (\ref{Notations-1-dbH}), there exists $\{\e_n\}_{n\geq1}$ such that
$$\lim\limits_{n\rightarrow\infty}\e_n=0,\ \ \lim\limits_{n\rightarrow\infty}\Big[\dbH(\e_n,X_1)-\dbH(\e_n,Y_1)\Big]=0.$$
\el

\begin{proof}
 \rm At first, for any $\e>0$ we make the following observation
\bel{X-1-H-h-1}\ba{ll}
\ns\ds \dbH(\e,X_1)=
\dbE\int_{\t+\e}^T \big[\frac 1 \e X_1(s)^{\top}\bar H_{xx}(s)X_1(s)
-Y_1(s)^{\top}\bar H_{xx}(s)Y_1(s)\big]ds\\
\ns\ds\qq\qq\ +\dbH(\e,Y_1)+\dbE\big[\frac 1 \e X_1(T)^{\top}\bar h_{xx}(T)X_1(T)
-Y_1(T)^{\top}\bar h_{xx}(T)Y_1(T)\big]\\
\ns\ds\qq\qq:=\dbH_1(\e)+\dbH(\e,Y_1)+\dbH_3(\e).
\ea\ee
We first treat $\dbH_1(\e)$. An easy calculation shows that
$$\ \ \ba{ll}
\ns\ds \Big|\dbH_1(\e)\Big| \leq \Big[ \dbE\int_{\t+\e}^T|\frac 1 \e X_1(s)X_1(s)^{\top}-Y_1(s)Y_1(s)^{\top}|^2ds\Big]^{\frac 1 2}
\Big[\dbE\int_{\t+\e}^T|\bar H_{xx}(s)|^2ds\Big]^{\frac 1 2}\\
\ns\ds\leq K\Big[\dbE\int_{\t+\e}^T|\e^{-\frac1 2}X_1(s)-Y_1(s)|^4ds\Big]^{\frac 1 4}
\Big[\dbE\int_{\t+\e}^T\big[|\e^{-\frac 1 2}X_1(s)|^4+|Y_1(s)|^4\big]ds\Big]^{\frac 1 4},
\ea$$
where we use the fact:
$$|a a^{\top}-b b^{\top}|=|[a-b]a^{\top}+b[a-b]^{\top}|\leq |a-b|[|a|+|b|],$$
with $a, b\in\dbR^n$. From (\ref{variational-equation-X-1}), (\ref{Y-1-e}) and (\ref{SVIE-estimate-2}), we immediately have
$$\ba{ll}
\ns\ds \sup_{t\in[\t+\e,T]}\dbE\big|\e^{-\frac 1 2}X_1(t)-Y_1(t)\big|^4\leq K\sup_{t\in[\t+\e,T]}\dbE|\e^{-\frac 1 2}\rho_1(t)-\varrho_1(t)|^4\\
\ns\ds\leq K\e^{-2}\dbE\Big[\int_\t^{\t+\e}|X_1(s)|^2ds\Big]^{2}+K\e^{-2}\sup_{t\in[\t+\e,T]}
\dbE\Big[\int_\t^{\t+\e}|\d\bar\si(t,s)-\d\bar\si(t,\t)|^2ds\Big]^{2}.
\ea$$
For the first term, denoted by $\cG_1(\e)$, on the right hand, we obtain $
\lim\limits_{\e\rightarrow0}\cG_1(\e)=0
$ by Lemma \ref{Lem-estimate-Appendix}.
As to the second term, $\cG_2(\e)$, thanks to Lemma \ref{Lemma-Banach-space-result}, there exists $\{\e_n\}_{n\geq1}$ such that
$$\ba{ll}
\ns\ds   \cG_2(\e_n)\leq \e_n^{-1}\dbE\int_\t^{\e_n+\t} \|\d\bar\si(\cd,s)-\d\bar\si(\cd,\t)\|_{C([0,T];\dbR^n)}^4ds\rightarrow0,\qq n\rightarrow\infty.
\ea$$
Consequently, for such $\{\e_n\}$, we conclude that
$$\lim_{n\rightarrow\infty}\sup\limits_{t\in[\t+\e_n,T]}\dbE\big|{\e_n}^{-\frac 1 2}X_1(t)-Y_1(t)\big|^4=0.$$
As a result, $\lim\limits_{n\rightarrow\infty}\Big| \dbH_1(\e_n)\Big|=0$.

Similarly we prove that $\lim\limits_{n\rightarrow\infty}\Big| \dbH_3(\e_n)\Big|=0$.
The conclusion is established via (\ref{X-1-H-h-1}).
\end{proof}

\mds

To treat $Y_1(\cd)$ in  Lemma \ref{Lemma-step-2-1}, for $\varrho_2(\cd):=\d\bar\si(\cd,\t)$, we need $Y_2(\cd)$ on $[t+\e,T]$,
\bel{Y-2-definition}\ba{ll}
\ns\ds Y_2(t)=\varrho_2(t)+\int_{\t+\e}^t\bar b_x(t,s)Y_2(s)ds+\int_{\t+\e}^t\bar\si_x(t,s)Y_2(s)dW(s).
\ea\ee
The solvability of $Y_2(\cd)\in C_{\dbF}([\t+\e,T];L^{4+\kappa}(\Omega;\dbR^n))$ is followed by the similar procedures as that of $Y_1(\cd)$ in (\ref{Y-1-e}).
Moreover, by the uniqueness in $C_{\dbF}([\t+\e,T];L^4(\Omega;\dbR))$, we see that
$$\left\{\ba{ll}
\ns\ds Y_1(\cd)=\e^{-\frac 1 2}(W(\t+\e)-W(\t))Y_2(\cd),\\
\ns\ds\dbH(\e,Y_1)=\dbE\Big[\e^{-1}|W(\t+\e)-W(\t)|^2\cd F_1^{\d\bar\si,\d\bar\si}(\t+\e)\Big],
\ea\right.$$
 where
$$\ba{ll}
\ns\ds  F_1^{\d\bar\si,\d\bar\si}(\t+\e):=\dbE_{\t+\e}\Big[\int_{\t+\e}^T Y_2(s)^{\top}\bar H_{xx}(s)Y_2(s)ds+ Y_2(T)^{\top}\bar h_{xx}(T)Y_2(T)\Big].
\ea$$
We establish the following result via $Y_2(\cd)$.

\mds

\bl\label{Lemma-step-2-2}
Given $\dbX(\cd)$, $Y_1(\cd)$, $\dbH(\e,Y_1)$ in (\ref{F-d-bar-si}), (\ref{Y-1-e}), (\ref{Notations-1-dbH}), we have
\bel{Y-1-Y-2-dbX}\ba{ll}
\ns\ds \lim_{\e\rightarrow0}\dbH(\e,Y_1) =\dbE\int_\t^T\dbX(s)^{\top}\bar H_{xx}(s)\dbX(s)ds+\dbE\big[\dbX(T)^{\top}\bar h_{xx}(T)\dbX(T)\big].
\ea\ee
\el
\begin{proof}
\rm Recalling $F^{\d\bar\si,\d\bar\si}(\cd)$ in (\ref{F-d-bar-si}), we shall derive the conclusion if
$$\ba{ll}
\ns\ds  \lim_{\e\rightarrow0}\dbE\Big[\e^{-1}|W(\t+\e)-W(\t)|^2 \big[F_1^{\d\bar\si,\d\bar\si}(\t+\e)
-F^{\d\bar\si,\d\bar\si}(\t)\big]\Big]=0.
\ea$$
To obtain this result, using similar ideas as in (\ref{Brownian-motion-formula-application}), we only prove
$$\lim_{\e\rightarrow0}\dbE\big|F_1^{\d\bar\si,\d\bar\si}(\t+\e)
-F^{\d\bar\si,\d\bar\si}(\t)\big|^{p_0}=0, \ \
p_0=\frac{2(4+\kappa)}{8+\kappa}\in(1,2).$$
To this end, by defining
$$\ba{ll}
\ns\ds  \Theta_1(\t+\e):= \int_{\t+\e}^T Y_2(s)^{\top}\bar H_{xx}(s) Y_2(s)ds+Y_2(T)^{\top}\bar h_{xx}(T)Y_2(T),\\
\ns\ds
\Theta_2(\t):= \int_{\t}^T \dbX(s)^{\top}\bar H_{xx}(s)\dbX(s)ds+ \dbX(T)^{\top}\bar h_{xx}(T)\dbX(T),
\ea$$
we deduce that
$$\ba{ll}
\ns\ds  \dbE\big|F_1^{\d\bar\si,\d\bar\si}(\t+\e)
-F^{\d\bar\si,\d\bar\si}(\t)\big|^{p_0}\equiv\dbE|\dbE_{\t+\e}\Th_1(\t+\e)-\dbE_\t\Th_2(\t)|^{p_0}\\
\ns\ds\leq K \dbE\big| \Th_1(\t+\e)-\Th_2(\t) \big|^{p_0}+K\dbE\big|\dbE_{\t+\e}\Th_2(\t)-\dbE_\t\Th_2(\t)\big|^{p_0}.
\ea$$
Hence it is suffice to prove the terms on right hand approach to zero as $\e\rightarrow 0$.

Notice that
\bel{estimate-Th-1-Th-2}\ba{ll}
%
%
\ns\ds \dbE\big| \Th_1(\t+\e)-\Th_2(\t) \big|^{p_0}\\
\ns\ds \leq K\dbE\Big\{\Big|\int_{\t+\e}^T|Y_2(s)Y_2(s)^{\top}-\dbX(s)\dbX(s)^{\top}|^2ds\Big|^{\frac {p_0} 2} \Big[\int_{\t+\e}^T|\bar H_{xx}(s)|^2ds\Big]^{\frac {p_0} 2}\Big\}\\
\ns\ds \q\ +K\dbE\Big\{\Big[\int_\t^{\t+\e}|\dbX(s)\dbX(s)^{\top}|^2ds\Big]^{\frac {p_0} 2}\Big[\int_\t^{\t+\e}|\bar H_{xx}(s)|^2ds\Big]^{\frac {p_0} 2}\Big\}\\
\ns\ds \q\ +K \dbE\Big\{|Y_2(T)Y_2(T)^{\top}-\dbX(T)\dbX(T)^{\top}|^{p_0}|\bar h_{xx}(T)|^{p_0} \Big\},
\ea\ee
where $x^{\top}A x=Tr(xx^{\top}A)$ with $x\in\dbR^n$ and $A\in \dbR^{n\times n}$.

As to the first term of (\ref{estimate-Th-1-Th-2}), denoted by $\cM_1(\e)$,
\bel{estimate-Th-1-Th-2-first-term} \ba{ll}
\ns\ds\cM_1(\e)\leq K\Big[\dbE\int_{\t+\e}^T|\bar H_{xx}(s)|^2ds\Big]^{\frac {p_0} 2} \Big[\dbE\Big(\int_{\t+\e}^T|Y_2(s)Y_2(s)^{\top}-\dbX(s)\dbX(s)^{\top}|^2ds\Big
)^{ p_0^* }\Big]^{\frac{2-{p_0}}{2}}\\
\ns\ds\qq\q\leq K\Big[\dbE\Big(\int_{\t+\e}^T\big[|Y_2(s)|+|\dbX(s)|\big]^4ds\Big)^{\frac { p_0 ^*}{2}}\Big(\int_{t+\e}^T|Y_2(s)-\dbX(s)|^4ds\Big)^{\frac { p_0 ^*}{2}}\Big]^{\frac{2-{p_0}}{2}}\\
\ns\ds\qq\q\leq K\Big[\dbE \int_{\t+\e}^T\big[|Y_2(s)|+|\dbX(s)|\big]^{4 p_0 ^*}ds \Big]^{\frac {2-{p_0}}{4}}
\Big[\dbE \int_{\t+\e}^T|Y_2(s)-\dbX(s)|^{4 p_0 ^*}ds \Big]^{\frac {2-{p_0}}{4}},
\ea\ee
where $ p_0 ^*:= \frac{p_0}{2-p_0}>1$.

 As to the second term of (\ref{estimate-Th-1-Th-2}), denoted by $\cM_2(\e)$,
\bel{Theta-1-Theta-2-middle}\ba{ll}
\ns\ds\cM_2(\e)\leq\Big[\dbE\int_\t^{\t+\e}|\bar H_{xx}(s)|^2ds\Big]^{\frac {p_0} 2}\cd \Big[\dbE\Big(\int_\t^{\t+\e}|\dbX(s)|^4ds\Big)^{ p_0 ^*}\Big]^{\frac {2-p}{2}}\rightarrow0,\ \ \e\rightarrow0.
\ea\ee

As to the third term of (\ref{estimate-Th-1-Th-2}), denoted by $\cM_3(\e)$,
\bel{estimate-Th-1-Th-2-third-term}\ba{ll}
\ns\ds\cM_3(\e)\!\leq\! K\!\Big[\!\dbE|\bar h_{xx}(T)|^2\!\Big]^{\frac {p_0} 2} \Big[\!\dbE\big[|Y_2(T)|+|\dbX(T)|\big]^{4 p_0^*} \dbE|Y_2(T)\!-\!\dbX(T)|^{4 p_0^*}\!\Big]^{\frac {2-p_0}{4}}.
\ea\ee
Since $4p_0^*=\frac{4p_0}{2-p_0}=4+\kappa$, we first estimate $\dbE|Y_2(\cd)-\dbX(\cd)|^{4+\kappa}$ which is similar as (\ref{SVIE-estimate-2}) or (\ref{Gronwall-SVIE-estimate}),
$$\ba{ll}
\ns\ds \sup_{r\in [\t+\e,T]}\dbE|Y_2(r)-\dbX(r)|^{4+\kappa}\\
\ns\ds\!\leq\! K\!\sup_{t\in[\t+\e,T]}\dbE\Big|\int_\t^{\t+\e}\bar b_{x}(t,s)\dbX(s)ds\!+\!\int_\t^{\t+\e}\bar\si_x(t,s)\dbX(s)dW(s)\Big|^{4+\kappa}\!\rightarrow\!0,\ \ \e\rightarrow0.
\ea$$
Therefore, from (\ref{estimate-Th-1-Th-2-first-term}), (\ref{estimate-Th-1-Th-2-third-term}), $\lim\limits_{\e\rightarrow0}\big[\cM_1(\e)+\cM_3(\e)\big]=0$.
Considering (\ref{Theta-1-Theta-2-middle}), one has
$$\lim\limits_{\e\rightarrow0}\dbE\big| \Th_1(\t+\e)-\Th_2(\t) \big|^{p_0}=0.$$

Our remaining aim is $\lim\limits_
{\e\rightarrow0}\dbE\big|\dbE_{\t+\e}\Th_2(\t)-\dbE_\t\Th_2(\t)\big|^{p_0}=0.$
By Lemma \ref{Lemma-SVIEs}
$$\ba{ll}
\ns\ds \dbE|\Th_2(\t)|^{p_0}\leq K\dbE\Big[\int_{\t}^T|\dbX(s)|^2|\bar H_{xx}(s)|ds\Big]^{p_0}+K\dbE\Big[|\dbX(T)|^{2{p_0}}|\bar h_{xx}(T)|^{p_0}\Big]\\
\ns\ds\qq\qq\q \leq K\Big[\dbE\int_\t^T|\d\bar\si(s,\t)|^{\frac{4{p_0}}{2-{p_0}}}ds\Big]^{\frac{2-{p_0}}{2}}+K \Big[\dbE |\d\bar\si(T,\t)|^{\frac{4p}{2-{p_0}}} \Big]^{\frac{2-{p_0}}{2}}<\infty.
\ea$$
Because $\dbE_{\t+\e}\Th_2(\t)\rightarrow \dbE_{\t}\Th_2(\t)$, a.s., $\e\rightarrow0$, and for any $r\in[\t,T]$,
$$\ba{ll}
\ns\ds \dbE|\dbE_r\Th_2(\t)|^{p_0}\leq\dbE\sup_{r\in[\t,T]}\dbE_r|\Th_2(\t)|^{p_0}\leq \frac{{p_0}}{{p_0}-1}\dbE|\Th_2(\t)|^{p_0}<\infty.
\ea$$
One has the desired conclusion by dominated convergence theorem.
\end{proof}

Now it is time for us to show the proof of Lemma \ref{Theorem-step-2}.

\begin{proof}
By above (\ref{variational-equation-X-1}), (\ref{Notations-1-dbH}), one has,
$$\ba{ll}
\ns\ds \frac 1 \e\Big[\dbE\int_0^T X_1(s)^{\top}\bar H_{xx}(s)X_1(s)ds+\dbE\big[X_1(T)^{\top}\bar h_{xx}(T)X_1(T)\big]\Big]\\
\ns\ds=\frac 1 \e \dbE\int_\t^{\t+\e} X_1(s)^{\top}\bar H_{xx}(s)X_1(s)ds+\dbH(\e,X_1).
\ea$$
For the first term on right hand, denoted by $\cQ(\e,X_1)$, from Lemma \ref{Lem-estimate-Appendix} we see that, $\lim\limits_{\e\rightarrow0}\cQ(\e,X_1)=0$.
We thus derive the conclusion by Lemma \ref{Lemma-step-2-1} and Lemma \ref{Lemma-step-2-2}.
\end{proof}

\subsection{Maximum principles of optimal control problems for SVIEs}

We present the first main result of this paper, the proof of which is based on the arguments from (\ref{variational-equations}) to (\ref{quadratic-form-X-1}), as well as Lemma \ref{SMP-step-1-main}, Lemma \ref{Theorem-step-2}.
\mds

\bt\label{SMP-main}
 Let (H1)-(H2) hold and $(\bar
X(\cd),\bar u(\cd))$ be an optimal pair. Then
$$\min\limits_{u\in U}\sH(t,u)=\sH(t,\bar u(t))=0,\ \ \dbP-\text{a.s.}, \ \ t\in[0,T],\ \ \text{a.e.} $$
where
\bel{maximum-1}\ba{ll}
\ns\ds \sH(t,u):= \D H^{\e}(t)+\frac 1 2 \d\bar\si(T,t)^{\top} [\cB_1(t)\D\bar\si(\cd,t)]\!+\!\frac 1 2 \int_0^T\d\bar\si(s,t)^{\top} \big[\cB_2(t)\D\bar\si(\cd,t)\big](s)ds,
\ea\ee
$\D H^{\e}(\cd),$ $\d\bar\si(T,t)$, $\D\bar\si(\cd,t)$ are in (\ref{bar-H-xx-definition}), (\ref{d-bar-si-1}),
  $\cB_1(\cd),$ $\cB_2(\cd)$
satisfy (\ref{boundedness-db-B-2-1-1-1}), (\ref{uniqueness-of-cB-1-2}).
\et

\mds

Above $\cB_1, \ \cB_2$ are called the $second$-$order$ $operator$-$valued$ $adjoint$ $processes$ of our optimal control problem.

If $\d\bar\si(\cd,\t)$ degenerates into $\cF_\t$-measurable random variable, using again the arguments from (\ref{variational-equations}) to (\ref{quadratic-form-X-1}), and Lemma \ref{Lemma-special-sigma}, Lemma \ref{Theorem-step-2}, we have the second main result in this article,

\mds

\bt\label{SMP-main-2}
 Let (H1)-(H2) hold with $\d\bar\si(\cd,\t)\equiv \d\bar\si(\t)$, and $(\bar
X(\cd),\bar u(\cd))$ be an optimal pair. Then
$$\min\limits_{u\in U}\sH_0(t,u)=\sH_0(t,\bar u(t))=0,\ \ \dbP-\text{a.s.},\ \ t\in[0,T],\ \ \text{a.e.}
$$
 where
\bel{maximum-1-SDE-special}\ba{ll}
\ns\ds \sH_0(t,u):=\D H^{\e}(t)+\frac 1 2 \d\bar\si(t)^{\top} \cB_3(t) \d\bar\si(t),
 \ \ u\in U,\
\ea\ee
$ \cB_3(\cd)$, $\D H^{\e}(\cd)$ are defined in Lemma \ref{SMP-step-1-SDE-main}, and (\ref{bar-H-xx-definition}), respectively.
\et

\mds

$\cB_3(\cd)$ is referred as the $\dbR^{n\times n}$-$valued$ $second$-$order$ $adjoint$ $process$ under this framework.

\mds

\begin{remark}\label{Special-sigma-suprising}
Suppose $\si_1,\ \si_2$ are two functions that satisfy the same requirements as $\si$ in (H1), and
\bel{Special-sigma-setting-2}\ba{ll}
\ns\ds
\si(t,s,x,u):=\si_1(t,s,x)+\si_2(s,x,u),\qq t,s\in[0,T],\ \ x\in\dbR^n,\ \ u\in U.
\ea\ee
For the previous $\d\bar\si(t,\t)$ in (\ref{d-bar-si-1}) with $ t,\ \t\in[0,T],$
\bel{d-bar-si-1-dbR}
\ba{ll}
\ns\ds
\d\bar\si(t,\t)\!:=\! \si(t,\t,\bar X(\t),u)\!-\!\si(t,\t,\bar X(\t),\bar u(\t))\!\equiv \! \d\bar\si(\t),\   u\in U.
\ea\ee
The corresponding maximum principle is easy to see in terms of Theorem \ref{SMP-main-2}.
\end{remark}

\mds

Next we discuss several special cases.

\subsubsection{State-independent diffusion and drift terms}

If both $b$ and $\si$ do not depend on $x$, then
$$\d\bar\si(t,\t):=\si(t,\t, u)-\si(t,\t,\bar u(\t)),\ \ t,\t\in[0,T],\ \ u\in U,$$
and (\ref{first-order-adjoint-equation}), (\ref{Hamiltonian-function}) change accordingly. In addition, for $u\in U$, $\bar H_{xx}(t,x,\bar X(\cd),\bar Y(\cd),\bar Z(\cd,t),u)$ is bounded.

We define two operator-valued processes $\cB'_1,\ \cB_2'$, i.e.  for $t\in[0,T]$, $s\in[0,T]$,
$$\ba{ll}
\ns\ds  \cB_1'(t )\D\bar\si(\cd,t):=\dbE_t[\bar h_{xx}(T)]\d\bar\si(T,t),\ \ \big[\cB_2'(t )\D\bar\si(\cd,t)\big](s):=\dbE_t[\bar H_{xx}(s)]\d\bar\si(s,t)I_{[t,T]}(s).
\ea$$
It is a direct calculation that $\cB_1'$ and $\cB_2'$ satisfy (\ref{boundedness-db-B-2-1-1-1}) and (\ref{main-result-1}). By the uniqueness in (\ref{uniqueness-of-cB-1-2}) and Theorem \ref{SMP-main}, for any $ u\in U$,
\bel{maximum-2}\ba{ll}
\ns\ds \D H^\e(t)\!+\!\frac 1 2 \d\bar\si(T,t)^{\top}\dbE_t[\bar h_{xx}(T)]\d\bar\si(T,t)\!+\!\frac 1 2\! \int_t^T\!\d\bar\si(s,t)^{\top}\dbE_t\bar H_{xx}(s)\d\bar\si(s,t)ds\!\geq\!0. \  a.s.
\ea\ee

\mds

\bc\label{State-independent-case}
Suppose $b$ and $\sigma$ do not rely on $x$, and $(\bar X(\cd),\bar u(\cd))$ is optimal pair.
Then for almost $t\in[0,T]$, (\ref{maximum-2}) holds true.
\ec

\mds

When the state equation is linear and the cost functional is quadratic, similar conclusion was obtained in \cite{Wang-2016-submitted}. In other words, our Corollary \ref{State-independent-case} extends theirs into the nonlinear setting.

\subsubsection{The convex control region}

Suppose (H1), (H2) hold with convex $U$. Moreover, the following maps are continuous differentiable,
\bel{convex-domain-differentiability}\ba{ll}
\ns\ds u\rightarrow \big(l(t,x,u),\ b(s,t,x,u),\ \si(s,t,x,u)\big),\ \ x\in\dbR^n,\ u\in U,\ \ s, t\in[0,T].
\ea\ee
For $t\in[0,T]$, a.e. $\o\in\Omega$, a.s., we define $v:=\bar u(t)+\e [u-\bar u(t)]$ with $u\in U$. The convexity of $U$ shows that $v\in U$. From (\ref{maximum-1}) we have,
$$\ba{ll}
\ns\ds 0\!\leq\!\frac{ \sH(t,v)}{\e}\!\leq \! \lan H_u(t,\bar X(t),\bar X(\cd),\bar Y(\cd),\bar Z(\cd,t),\bar u(t)+\theta\e(u-\bar u(t))), u-\bar u(t) \ran \!+\!K(t)\e,
\ea$$
where $K(\cd)$ is a process, $0<\theta<1$ and $H_u$ is the partial derivative with respect to $u$.
Let $\e\rightarrow0$, one has
\bel{maximum-1-convex}\ba{ll}
\ns\ds \lan H_u(t,\bar X(t),\bar X(\cd),\bar Y(\cd),\bar Z(\cd,t),\bar u(t)), u-\bar u(t)\ran\geq0.
\ \ \ a.s.
\ea\ee
\begin{corollary}
Let (H1), (H2), (\ref{convex-domain-differentiability}) hold true and $(\bar
X(\cd),\bar u(\cd))$ be optimal with convex $U$. Then there exists a pair of $(\bar Y(\cd),\bar Z(\cd,\cd))$ satisfying
(\ref{first-order-adjoint-equation}) such that for almost $t\in[0,T]$, $u\in U,$ (\ref{maximum-1-convex}) is satisfied.
\end{corollary}

\mds

Above convex region case was studied in e.g. \cite{Yong 2006}, \cite{Yong 2008}. Therefore, our study (i.e. Theorem \ref{SMP-main}) extends theirs into the non-convex setting.

\subsubsection{The linear quadratic case.}

We discuss the linear quadratic optimal control problem where $U:=\dbR^m$,
$$\ba{ll}
\ns\ds b=[A_1(t,s)x+B_1(t,s)u],\ \ \si=[A_2(t,s)x+B_2(t,s)u],\ \\
\ns\ds  h=\frac 1 2 x^{\top}G x, \ \ l=\frac 1 2\big[x^{\top}Qx+2u^{\top}Sx+u^{\top}Ru\big].
\ea$$
Here $A_1,B_1, A_2, B_2, Q, S, R, G$ are bounded and random such that for some modulus function $\rho(\cd)$,
$$
 |f(t,s)-f(t',s)| \leq\rho(|t-t'|),\ \ t,t',s\in[0,T],\ \ f:= A_1,B_1, A_2, B_2.
$$
For optimal $(\bar X,\bar u)$, we define
\bel{Linear-quadratic-Hamiltion-function-definition}\ba{ll}
\ns\ds  \sH_1(t):=  S(t)\bar X(t)+B_1(T,t)^{\top}\dbE_t\big[G\bar X(T)\big]+\dbE_t\int_t^TB_1(s,t)^{\top}\bar Y(s)ds  \\
\ns\ds\qq\qq +B_2(T,t)^{\top}\bar\pi(t)+\dbE_t\int_t^TB_2(s,t)^{\top}\bar Z(s,t)ds,
\ea\ee
where $(\bar Y,\bar Z,\bar\pi )$ is in (\ref{first-order-adjoint-equation}) accordingly.
Next we consider two special cases.

\mds

\bf Case I: \rm

Suppose $A_1\equiv A_2 \equiv0$. For any $u\in\dbR^m$, $t\in[0,T],$ a.e., (\ref{maximum-2}) becomes
$$\ba{ll}
\ns\ds
 [u-\bar u(t)]^{\top}\!\sH_1(t) \!+\!\frac{1}{2}\big[u^{\top}R(t)u\!-\!\bar u(t)^{\top}R(t)\bar u(t)\big]\!+\!\frac 1 2 [u\!-\!\bar u(t)]^{\top}\sB(t)[u-\bar u(t)]\geq0,\ a.s.
 \ea
$$
where $\sH_1(t)$ is in (\ref{Linear-quadratic-Hamiltion-function-definition}), $\sB(t):=\sB_1(t)+\sB_2(t)$,
$$\sB_1(t):= B_2(T,t)^{\top}\big[\dbE_tG\big] B_2(T,t),\ \ \sB_2(t):=\dbE_t\int_t^TB_2(s,t)^{\top}Q(s) B_2(s,t)ds.$$
%
Here $\sB_i$ are $\dbF$-adapted processes. By the arbitrariness of $u$, one has
$$\sH_1(t)+R(t)\bar u(t)=0,\ \ R(t)+\sB(t)\geq0,\ \ t\in[0,T]. \ \ a.s. \ \ a.e.$$
Notice that these two condition were also obtained in \cite{Wang-2016-submitted} with distinct approach.
%
%

\mds

\begin{remark}\label{Some-essential-ideas}
If for $t\in[0,T]$, $B_i(\cd,t)\equiv B_i(t)$, $i=1,2$. Then we have
$$\sB(t)=B_2(t)^{\top}P_2(t)B_2(t),\ \  P_2(t):=\dbE_t\Big[G+\int_t^TQ(s)ds\Big].$$
Here $P_2$, which is called second-order adjoint process, just satisfies (\ref{Second-order-adjoint-equations-SDEs}), see \cite{Wang-2016-submitted}. For this specification, we reveal two interesting facts.

It is our believe that $G$ and $Q(\cd)$ should originally play its own peculiar role in optimal control problem of SVIEs, such as above introduced processes $\sB_i$. In particular SDEs case, such independence disappears and their roles happen to be merged together in suitable manner, like above $P_2$. In other words, compared with the SDEs scenario with only one adjoint process $P_2$, a pair of different forms of adjoint processes are indeed required in SVIEs setting.

For optimal control problems of SDEs, one can directly introduce second-order adjoint equation (\ref{Second-order-adjoint-equations-SDEs}). However, as to the case of SVIEs, this adjoint equation idea does not work any more (see \cite{Wang-2016-submitted}). Actually, in our opinion, there is no way to construct the analogue version of (\ref{Second-order-adjoint-equations-SDEs}) here. One has to seek more fundamental, appropriate notion to get around this difficulty. In a nutshell, we need to use new stochastic processes to replace classical stochastic equations.

\end{remark}

\mds

\bf Case II: \rm

Suppose $B_2(t,\cd)\equiv B_2(\cd)$. Then $\d\bar\si( t)=B_2( t)[u-\bar u(t)]$, $t\in[0,T]$.
Recall Theorem \ref{SMP-main-2}, we have
$$\ba{ll}
\ns\ds
 [u-\bar u(t)]^{\top}\sH_1(t) +\frac{1}{2}\big[u^{\top}R(t)u-\bar u(t)^{\top}R(t)\bar u(t)\big]\\
\ns\ds\qq +\frac 1 2 [u-\bar u(t)]^{\top}B_2(t)^{\top}\sP(t)B_2(t)[u-\bar u(t)]\geq0,\
 \ea
$$
where $u\in\dbR^m,$ $\sH_1(\cd)$ is in (\ref{Linear-quadratic-Hamiltion-function-definition}), $\sP(\cd):=\big\{\sP_{ij}(\cd)\big\}_{n\times n}$,
$$\left\{\ba{ll}
\ns\ds \sP_{ij}(t)=\frac{1}{2}\dbE_t\Big[\sX^{e_j}(T)^{\top}G\sX^{e_i}(T)\Big]+\frac 1 2 \dbE_t\int_t^T \sX^{e_j}(s)^{\top}Q(s)\sX^{e_i}(s)ds,\ \ t\in[0,T],\\
\ns\ds \sX^{e_i}(r)=e_i+\int_t^rA_1(t,s)\sX^{e_i}(s)ds+\int_t^rA_2(t,s)\sX^{e_i}(s)dW(s),\ \ r\in[t,T].
\ea\right.$$
By the arbitrariness of $u\in\dbR^m$, one then obtains the following maximum condition
\bel{optimality-necessary-condition-SLQ-SVIEs}\ba{ll}
\sH_1(t)+R(t)\bar u(t)=0,\ \ R(t)+B_2(t)^{\top}\sP(t)B_2(t)\geq0,\ \ t\in[0,T].\ \ a.e.
\ea\ee
\br
If $A_1,\ A_2,\ B_1$ are independent of $t$, then according to Lemma \ref{Lemma-special-sigma} and the arguments in Subsection 3.1, $\sP(\cd)$ is just the unique solution of second-order adjoint equation (\ref{Second-order-adjoint-equations-SDEs}), and (\ref{optimality-necessary-condition-SLQ-SVIEs}) naturally reduce into the counterpart in SDEs case.
\er

\section{Concluding remarks}

This article is devoted to maximum principles of optimal control problems for SVIEs when the control region is arbitrary subset of $\dbR^m$ and diffusion relies on control variable. Some novelties are summed up as follows.

$\bullet$ We introduce a class of quadratic functionals associated with linear SVIEs, and represent them in two different ways. To our best, these conclusions are new and may have independent interests. When they are applied in optimal control problems of SDEs, the maximum principles can be established without It\^{o} formula and second-order adjoint equations.

%

$\bullet$ We establish two maximum principles in terms of operator-valued, and matrix-valued second-order adjoint processes, respectively.
For optimal control problems of SVIEs (\ref{FSVIE-3-1}), there is no existing paper treating the case of closed $U$ and state-dependent diffusion.
Nevertheless, it is one particular case of our study. Moreover, our conclusions can fully cover the SDEs case.

$\bullet$ We obtain some convincing arguments to show that the second-order adjoint equation idea actually fails in our SVIEs setting. Therefore, we propose appropriate second-order adjoint processes instead. In addition, unlike the classical scenario with only one second-order adjoint equation, here we have to rely on two second-order adjoint processes in the maximum conditions, which of course merge into the solution of second-order adjoint equation in particular SDEs setting.

We emphasize that Theorem \ref{SMP-main-2} is a refinement of Theorem \ref{SMP-main}. As a trade-off, some requirements are imposed (Remark \ref{Special-sigma-suprising}). Consequently, it still remains its importance to replace the pair of operator-valued processes in Theorem \ref{SMP-main} by explicit $\dbR^{n\times n}$-valued counterparts. We hope to discuss this topic in our forthcoming papers.

\section{Appendix}

\mds

To prove the existence of operator-valued processes in Lemma \ref{Lemma-existence-operator-valued}, we make some preparations in the sequel. The first two results are more or less standard in functional analysis.

\bl\label{Appendix-lemma-3}
Suppose $f:\dbB\mapsto\dbR$ is a bounded linear functional. Then there exists a unique $y(\cd)=(y_1,y_2(\cd))\in \dbB'$ such that
$$\left\{\ba{ll}
\ns\ds f(x)=\int_0^T x_2(t)^{\top}y_2(t)dt+x_1^{\top} y_1,\qq \forall x=(x_1,x_2(\cd))\in \dbB, \\
\ns\ds \|f\|_{\cL(\dbB;\dbR)}=\|y\|_{\dbB'}=\max\Big\{|y_1|,\Big[\int_0^T|y_2(t)|^{\frac 4 3}dt\Big]^{\frac 3 4}\Big\}.
\ea\right.$$
\el

\bl\label{Appendix-lemma-4}
Given Banach space $X$ with numerable dense subset $X_0$, suppose
$$X_1:=\big\{\sum\limits_{i=1}^n a_ix_i,\ a_i\in\dbQ,\ x_i\in X_0,\ n\in\dbN\ \big\},$$
with $\dbQ$ the set of rational number, $f:X_1\mapsto Y$ is map to Banach space $Y$ such that for constant $M>0$,
$$\ba{ll}
\ns\ds \|f(x)\|_{Y}\leq M\|x\|_{X}, \ \ \forall x\in X_1,\ \ f(ax+by)=a f(x)+b f(y),\ \  \forall a,b\in\dbQ.
\ea$$
Then there exists a unique bounded linear operator $F:X\mapsto Y$ satisfying
$$F(x)=f(x),\ \ x\in X_1,\ \ \|F\|_{\cL(X,Y)}\leq M.$$
\el

\bl\label{Lemma-existence-operator-process-Appendix}
Given constant $M>0$, suppose $f:\dbB_1\times \dbB_1\mapsto\dbR$ satisfies
\bel{fake-bilinear-1}\ba{ll}
\ns\ds \|f(x,y)\| \!\leq\! M\|x\|_{\dbB_1}\|y\|_{\dbB_1},\  \forall x,\ y\in \dbB_1, \  f(a\bar x+b\wt x,y)\!=\!a f(\bar x,y)\!+\!b f(\wt x,y), \\
\ns\ds  f(x,a\bar y\!+\!b\wt y)\!=\!a f(x,\bar y)\!+\!b f(x,\wt y),\  a,b\in\dbQ,\  x,\bar x, \wt x, y,\bar y, \wt y\in\dbB_1.
\ea\ee
Then there exist a unique pair of linear operators $\h\cB_1:\dbB\mapsto\dbR^n$, $\h\cB_2:\dbB\mapsto L^{\frac 4 3}(0,T;\dbR^n)$ such that
$$\left\{\ba{ll}
\ns\ds \|\cB_1\|_{\cL(\dbB;\dbR^n)}\leq M,\ \ \|\cB_2\|_{\cL(\dbB;L^{\frac 4 3}(0,T;\dbR^n))}\leq M,\\
\ns\ds f(x,y)=\int_0^Tx_2(t)^{\top}\big[\h\cB_2 y\big](t)dt+x_1^{\top}\big[\h\cB_1 y\big],\ \ \forall x,\ y\in \dbB_1.
\ea\right.$$
\el

\begin{proof}
\rm
$Step$ $1:$ Fix $y\in \dbB_1$, we define $\varphi_y:\dbB_1\mapsto\dbR$, i.e. $\varphi_y(x)=f(x,y)$, $x\in \dbB_1$. It is easy to see
$$\varphi_y(k\bar x+l\wt x)=k\varphi_y(\bar x)+l\varphi(\wt x),\ \ |\varphi_y(x)|\leq M\|x\|_{\dbB}\|y\|_{\dbB},\ \ k,l\in\dbQ, \ \ \bar x,\wt x\in \dbB_1.$$
According to Lemma \ref{Appendix-lemma-4}, there exists a unique bounded linear functional $\h \varphi_y: \dbB\mapsto\dbR$ such that
$$\h\varphi_y(x)=\varphi_y(x),\ \ \|\h\varphi_y\|_{\cL(\dbB;\dbR)}\leq \cM_y:= M\|y\|_{\dbB}, \ \ \forall x\in \dbB_1.$$
Consequently, on account of Lemma \ref{Appendix-lemma-3}, there exists a unique $y^*(\cd)=( y_1^*, y_2^*(\cd))\in \dbB'$ such that
$$\left\{\ba{ll}
\ns\ds \h\varphi_y(x)\!=\!x_1^{\top} y^*_1\!+\!\!\int_0^T x_2(t) ^{\top}y^*_2(t)dt,\ \ x=(x_1,x_2(\cd))\in \dbB,\\
\ns\ds\|\h\varphi_y\|_{\cL(\dbB;\dbR)}\!=\!\max \!\Big\{|  y^*_1|,\Big[\int_0^T| y^*_2(t)|^{\frac 4 3}dt\Big]^{\frac 3 4}\Big\}.
\ea\right.$$
$Step$ $2:$
We introduce operator $\cB:\dbB_1\mapsto \dbR^n\times L^{\frac 4 3}(0,T;\dbR^n)$, i.e.,
$$\cB y=\big([\cB y]_1,[\cB y]_2(\cd)\big),\ \ [\cB y]_1= y^*_1, \  \[\cB y]_2(\cd)= y^*_2(\cd), \ \ y\in\dbB_1.$$
The well-posedness of $\cB$ is obvious. We claim that $\cB$ is linear in following sense,
$$\cB(k\bar h+l \wt h)=k\cB\bar h+l\cB\wt h, \  \ k,\ l\in\dbQ, \  \ \bar h,\ \wt h\in \dbB.
$$

Notice that $[k\bar h+l\wt h]\in\dbB_1$. From $Step$ $1$ there exists a unique linear bounded $\h\varphi_{k\bar h+l\wt h}(\cd):\dbB\mapsto\dbR$ such that $\h\varphi_{k\bar h+l\wt h} =\varphi_{k\bar h+l\wt h} $ in $\dbB_1$. Since $\varphi_{k\bar h+l\wt h} \!=\!k\varphi_{\bar h} +l\varphi_{\wt h} $ in $\dbB_1$,
it then follows that $k\h\varphi_{\bar h}+l\h\varphi_{\wt h}$ is another linear bounded functional on $\dbB$ such that $\big[k\h\varphi_{\bar h}+l\h\varphi_{\wt h}\big] =\varphi_{k\bar h+l\wt h} $ in $\dbB_1.$
By virtue of the uniqueness in Lemma \ref{Appendix-lemma-4}, $ \h\varphi_{k\bar h+l\wt h} =\big[k\h\varphi_{\bar h}+l\h\varphi_{\wt h}\big].$

From Lemma \ref{Appendix-lemma-3} we see that there exists
$$\big[k\bar h+l\wt h\big]^*=\Big(\big[ k\bar h+l\wt h \big]^*_1,\big[ k\bar h+l\wt h \big]^*_2(\cd)\Big)\in \dbR^n\times L^{\frac 4 3}(0,T;\dbR^n)$$
such that,
$$\ba{ll}
\ns\ds  \h\varphi_{k\bar h+l\wt h}(x)\!=\!\Big[\big[ k\bar h\!+\!l\wt h \big]^*_1\Big]^{\top} x_1\!+\!\int_0^T\Big[\big[ k\bar h\!+\!l\wt h \big]^*_2(s)\Big]^{\top}x_2(s)ds,\ \ \forall x=(x_1,x_2(\cd))\in \dbB.
\ea$$
Using again Lemma \ref{Appendix-lemma-3},
$$\ba{ll}
\ns\ds  k\h\varphi_{\bar h}(x)+l\h\varphi_{\wt h}(x)=\Big[k\int_0^Tx_2(t)^{\top}\bar y^*_2(t)dt+k x_1^{\top} \bar y^*_1\Big]+\Big[l\int_0^Tx_2(t)^{\top}\wt y^*_2(t)dt+l x_1^{\top} \wt y^*_1\Big]\\
\ns\ds\qq\qq=\int_0^Tx_2(t)^{\top}\big[k\bar y^*_2(t)+l\wt y^*_2(t)\big]dt+ x_1^{\top} \big[k\bar y^*_1+l\wt y^*_1\big].
\ea$$
Therefore, by the arbitrariness of $x=(x_1,x_2(\cd))$,
$$\big[ k\bar h+l\wt h \big]^*_2 =k\bar y^*_2 +l\wt y^*_2,\ \ \big[ k\bar h+l\wt h \big]^*_1=k\bar y^*_1+l\wt y^*_1.$$
This directly leads to the desirable linearity of $\cB$.
%

$Step$ $3:$ We extend the definition of $\cB$ into $\dbB$. For any $y\in \dbB_1,$ from $Step$ $1,2$,
$$\|\cB y\|_{\dbB'}=\|\h\varphi_y\|_{\cL(\dbB;\dbR)}\leq M\|y\|_{\dbB}.$$
Hence from the linearity of $\cB$ and Lemma \ref{Appendix-lemma-4}, there exists a unique linear $\h\cB:\dbB\mapsto\dbB'$ such that
$$\h\cB(x)=\cB(x),\ \ x\in \dbB_1,\ \ \|\h\cB\|_{\cL(\dbB,L^{\frac 4 3}(0,T;\dbR^n))} \leq M.$$
Therefore, for any $x,\ y\in \dbB_1$,
\bel{f-x-y-varphi}\ba{ll}
\ns\ds f(x,y)=\varphi_y(x)=\h\varphi_y(x)=x_1^{\top}\big(\h\cB y\big)_1+\int_0^Tx_2(t)^{\top}\big(\h\cB y\big)_2(t)dt.
\ea\ee
$Step$ $4:$ To see the conclusions of $\h\cB_1$, $\h\cB_2$, for $y\in \dbB_1$, we define $\cB_1:\dbB_1\mapsto\dbR^n$, $\cB_2:\dbB_1\mapsto L^{\frac 4 3}(0,T;\dbR^n)$ as
$$\cB_1y=\big(\cB y\big)_1,\ \ (\cB_2y)(\cd)=\Big[\big(\cB y\big)_2\Big](\cd).$$
Of course, both of them are well-defined.

We look at their linearity with $k, l\in\dbQ$, $\bar h,\wt h\in \dbB_1$. By the linearity of $\cB$, we obtain
$$\cB_i(k\bar h+l\wt h)=\big[\cB(k\bar h+l\wt h)\big]_i=k\big(\cB\bar h\big)_i+l\big(\cB\wt h\big)_i,\ \ i=1,2.$$
On the other hand, for any $h\in \dbB_1$ we know that $|\cB_i h|\leq\|\cB h\|\leq M\|h\|_{\dbB},$ $i=1,2$.
By using Lemma \ref{Appendix-lemma-4} again, there exist a unique linear bounded $\h\cB_1:\dbB\mapsto\dbR^n$ and a unique linear bounded $\h\cB_2:\dbB\mapsto L^{\frac 4 3}(0,T;\dbR^n)$ such that $\h\cB_1 y=\cB_1 y$, $\h\cB_2 y=\cB_2 y$ with $y\in \dbB_1.$ As a result, from (\ref{f-x-y-varphi}) we obtain the conclusion.
\end{proof}

Using almost the same ideas as above, one obtain the following which is useful in proving Lemma \ref{Lemma-special-sigma}.

\bl\label{Lemma-existence-real-valued-process-Appendix}
Given positive constant $M$, suppose map $f:\dbQ^n\times \dbQ^n\mapsto\dbR$ satisfies
\bel{fake-bilinear-1-finite-dimensional}\ba{ll}
\ns\ds |f(x,y)|\! \leq\! M|x|_{\dbQ^n}|y|_{\dbQ^n},\  \forall x,y\in \dbQ^n, \  f(a\bar x+b\wt x,y)\!=\!a f(\bar x,y)+b f(\wt x,y), \\
\ns\ds  f(x,a\bar y+b\wt y)\!=\!a f(x,\bar y)\!+\!b f(x,\wt y),\  a,b\in\dbQ,\  x,\ \bar x,\ \wt x, \ y,\ \bar y,\ \wt y\in\dbQ^n,
\ea\ee
where $\dbQ^n$ is the set of $n$-dimensional rational vectors. Then there exist a unique $\dbR^{n\times n}$-valued matrix $\h\cB_3$ such that $|\h\cB_3|\leq M$, and
%
$
f(x,y)= x^{\top}\h\cB_3 y,$ $  x,y \in\dbQ^n.$

\el

\textbf{Acknowledgements.} The author highly appreciates the anonymous referees' constructive comments. He also gratefully acknowledges Professor Jiongmin Yong and Professor Xu Zhang for their valuable suggestions and helpful discussions.

\mds


\begin{thebibliography}{9}




\bibitem{Agram-Oksendal 2015}
{\sc N. Agram and B. {\O}ksendal},
{\em Mallivain calculus and optimal control of stochastic Volterra equations},
J. Optim. Theory Appl., (2015), DOI 10.1007/s10957-015-0753-5.


\bibitem{K. L Arrow 1964}
{\sc K. Arrow,} {\em Optimal capital policy, the cost of capital and myopic decision rules}, Ann. Inst. Stat. Math., 16 (1964), pp. 21--30.


\bibitem{Bakke 1974} {\sc V. Bakke}, {\em A maximum principle for an optimal control problem: with integral constraints}, J. Optim. Theory Appl. 13 (1974), pp. 32--55.



\bibitem{Bonaccorsi-Confortola-Mastrogiacomo 2012}
{\sc S. Bonaccorsi, F. Confortola and E. Mastrogiacomo}, {\em Stochastic control for stochastic Volterra equations with complete monotone kernels}, SIAM J. Control Optim., 50 (2012), pp. 748--789.


\bibitem{Bonnans et al 2013} {\sc J. Bonnans, X. Dupuis and C. De la Vega}, {\em First and second order optimality conditions for optimal control problems of state constrained integral equations}, J. Optim. Theory Appl. 159 (2013), pp. 1--40.

\bibitem{Cochran et al 1995} {\sc W. Cochran, J. Lee and J. Potthoff}, {\em Stochastic
Volterra equations with singular kernels}, \rm Stochastic Process Appl. 56 (1995), pp. 337--349.


\bibitem{De-Hoog-Weiss 1973} {\sc F. De Hoog and R. Weiss}, {\em On the solution of a Volterra integral
equation with a weakly singular kernel}, \rm SIAM J. Math. Anal. 4 (1973), pp. 561--573.

\bibitem{Dmitruk 2014} {\sc A. Dmitruk and N. Osmolovski}, {\em Necessary conditions for a weak minimum in optimal control problems with integral equations subject to state and mixed constraints}, \rm SIAM J. Control Optim. 52 (2014), pp. 3437--3462.

\bibitem{Du-Meng 2013} {\sc K. Du and Q. Meng}, {\em A maximum principle for optimal control of stochastic evolution equations}, \rm SIAM J. Control Optim. 51 (2013), pp. 4343--4362.






\bibitem{Friedman 1964} {\sc A. Friedman}, {\em Optimal control for hereditary processes}, \rm Arch. Rat. Mech. Anal. 15 (1964), pp. 396--416.

\bibitem{Fuhrman et al 2013} {\sc M. Fuhrman, Y. Hu and G. Tessitore}, {\em Stochastic maximum principle for optimal control of SPDEs}, \rm Appl. Math. Optim. 68 (2013), pp. 181--217.

\bibitem{Halanay 1968} {\sc A. Halanay}, {\em Optimal control for systems with time-lag}, \rm SIAM J. Control 6 (1968), pp. 215--234.



\bibitem{Hartl 1984}
{\sc R. Hartl}, {\em Optimal dynamic advertising policies for hereditrary processes}, J. Optim. Theory Appl. 43 (1984), pp. 51--72.

\bibitem{Hritonenko-Yatsenko 2008} {\sc N. Hritonenko and Y. Yatsenko}, {\em Optimal control of Solow vintage capital model with nonlinear utility}, \rm Optimization, 57 (2008), pp. 581--592.


\bibitem{S.C. Kou 2008} {\sc S. Kou}, {\em Stochastic modeling in nanoscale biophysics: subdiffusion within proteins}, \rm Ann. Appl. Stat. 2 (2008), pp.
501--535.



\bibitem{Lin-Yong-2018} {\sc P. Lin and J. Yong}, {\em Controlled singular Volterra integral equations
and Pontryagin maximum principle}, arXiv:1712.05911v1.


\bibitem{Lv-Zhang 2015} {\sc Q. L\"{u} and X. Zhang}, {\em General Pontryagin-type stochastic maximum principle and backward stochastic evolution equation in infinite dimensions}, Springer Briefs in Mathematics, 2014.

\bibitem{Kamien-Muller 1976} {\sc M. Kamien and E. Muller}, {\em Optimal control with integral
state equations}, \rm Rev. Econ. Stud., \rm 43  (1976), pp.
469--473.




\bibitem{Oksendal-Zhang 2010} {\sc B. {\O}ksendal and T. Zhang}, {\em Optimal control with partial information for stochastic Volterra equations}, \rm Int. J. Stoch. Anal., (2010), doi:10.1115/2010/329185



\bibitem{Pardoux-Protter 1990}
{\sc E. Pardoux and P. Protter}, {\em Stochastic Volterra equations
with anticipating coefficients}, Ann. Probab., 18 (1990), pp. 1635--1655.



\bibitem{Peng90}
{\sc S. Peng}, {\em A general stochastic maximum principle for optimal control problems}, SIAM J. Control Optim., 28 (1990), pp. 966--979.


%


%




\bibitem{Shi-Wang-Yong 2015}
{\sc Y. Shi, T. Wang and J. Yong}, {\em Optimal control problems of forward-backward
stochastic Volterra integral equations}, Math. Control Relat. Fields, 5 (2015), pp.  613--649.


\bibitem{Vinokurov 1969} {\sc V. Vinokurov}, {\em Optimal control of processes described by integral equations}, I, II, III, \rm Izv. Vys\v{s}. U\v{c}ebn. Zaved. Matematika 7, 21--33; 8, 16--23; 9, 16--25; (in Russian) English transl. in SIAM J. Control 7 (1967), pp. 324--336, 337--345, 346--355.



\bibitem{Wang-2016-submitted} {\sc T. Wang}, {\em Linear quadratic control problems of stochastic Voltera integral equations}, to appear in ESAIM: Control Optim. Cal. Var. DOI: https://doi.org/10.1051/cocv/2017002.  \rm





\bibitem{Wang-Yong 2013}
{\sc T. Wang and J. Yong}, {\em Comparison theorems for backward
stochastic Volterra integral equations}, Stochastic Process Appl., 125 (2015), pp. 1756--1798.


\bibitem{Wang-Zhang 2016} {\sc T. Wang and H. Zhang}, {\em Optimal control problems of forward-backward stochastic Volterra integral equations with closed control regions}, \rm SIAM J. Control Optim. 55 (2017), pp. 2574--2602.



\bibitem{Yong 2006} {\sc J. Yong}, {\em Backward stochastic Volterra integral
equations and some related problems}, Stoch. Process Appl.,
 116 (2006), pp. 779--795.
%


\bibitem{Yong 2008}
{\sc J. Yong}, {\em Well-posedness and regularity of
backward stochastic Volterra integral equation}, Probab. Theory
Relat. Fields, 142 (2008), pp. 21--77.




\bibitem{Yong-Zhou 1999}
{\sc J. Yong and X. Zhou}, {\em Stochastic Controls: Hamiltonian Systems and HJB
Equations}, Springer-Verlag, New York, Berlin, 2000.











\end{thebibliography}
\end{document}